\documentclass[11pt]{amsart}
\usepackage{amssymb,mathrsfs,stmaryrd}	
\usepackage[all]{xy}
\usepackage{graphicx}
\usepackage[mathcal]{euscript}
\usepackage{verbatim}
\usepackage{tikz}
     \usetikzlibrary{positioning}
     \usepackage{subcaption}

 \usepackage{bm}
\usepackage[inline]{enumitem}

\usepackage[colorlinks,citecolor=blue]{hyperref}

\renewcommand{\and}{\qquad\text{and}\qquad}

\let\oldmarginpar\marginpar
\renewcommand\marginpar[1]{\-\oldmarginpar[\raggedleft\footnotesize #1]%
{\raggedright\footnotesize #1}}

\usepackage{mathtools}
\usepackage{manfnt}

\usepackage{color}

\definecolor{jade}{rgb}{0.10, 0.56, 0.42}
\definecolor{cerise}{rgb}{0.87, 0.19, 0.39}
\usetikzlibrary{arrows.meta}
\tikzset{>={latex[width=3mm,length=3mm]}}

\usetikzlibrary{arrows,backgrounds, fit, calc, positioning}
\usepackage{multirow}

\theoremstyle{definition}
\newtheorem*{thm*}{Theorem}
\newtheorem{thm}{Theorem}[section]
\newtheorem{definition}[thm]{Definition}
\newtheorem{conjecture}[thm]{Conjecture}
\newtheorem{lem}[thm]{Lemma}
\newtheorem{prop}[thm]{Proposition}
\newtheorem{corollary}[thm]{Corollary}

\newtheorem*{Three}{Three Questions}
\newtheorem*{thm:mainthm}{Theorem \ref{thm:mainthm}}
\newtheorem*{thm:rigidthm}{Theorem \ref{thm:rigidthm}}
\theoremstyle{remark}
\newtheorem{rmk}[thm]{Remark}
\newtheorem{example}[thm]{Example}

\numberwithin{equation}{section}

\newcommand{\N}{\mathbb{N}}
\newcommand{\Z}{\mathbb{Z}}
\newcommand{\Q}{\mathbb{Q}}
\newcommand{\R}{\mathbb{R}}
\newcommand{\C}{\mathbb{C}}

\newcommand{\pp}{\mathbb{P}^1}
\newcommand{\End}{\mathrm{End}}
\newcommand{\Ad}{\mathrm{Ad}}
\newcommand{\ad}{\mathrm{ad}}

\newcommand{\sC}{\mathscr{C}}
\newcommand{\sM}{\mathscr{M}}
\newcommand{\orb}{\mathscr{O}}
\newcommand{\borb}{\bm{\cO}}

\newcommand{\cO}{\mathcal{O}}

\newcommand{\cA}{\mathcal{A}}
\newcommand{\cB}{\mathcal{B}}
\newcommand{\cC}{\mathcal{C}}
\newcommand{\cL}{\mathcal{L}}
\newcommand{\cE}{\mathcal{E}}

\newcommand{\ftype}{\mathcal{D}}
\newcommand{\ftyperam}{\mathcal{F}}
\newcommand{\bftype}{\bm{\mathcal{D}}}
\newcommand{\bftyperam}{\bm{\mathcal{F}}}

\newcommand{\DS}{\mathrm{DS}}

\newcommand{\bfF}{\mathbf{F}}
\newcommand{\bfD}{\mathbf{D}}
\newcommand{\bfO}{\mathbf{O}}
\DeclareMathOperator{\M}{\mathcal{M}}

\DeclareMathOperator{\fg}{\mathfrak{g}}
\DeclareMathOperator{\fgl}{\mathfrak{gl}}

\DeclareMathOperator{\GL}{GL}
\DeclareMathOperator{\SL}{SL}
\DeclareMathOperator{\Sp}{Sp}
\DeclareMathOperator{\Irr}{Irr}
\DeclareMathOperator{\Rep}{Rep}

\newcommand{\Stor}[1]{S_{{#1}}}
\newcommand{\symmetric}[1]{\mathfrak{S}_{{#1}}}

\newcommand{\Iwa}{I}

\newcommand{\Galgp}{\mathcal{I}}

\newcommand{\Para}{P}

\newcommand{\ft}{\mathfrak{t}}

\newcommand{\fcox}{\mathfrak{c}}
\newcommand{\diag}{\mathrm{diag}}
\newcommand{\fb}{\mathfrak{b}}
\newcommand{\fu}{\mathfrak{u}}

\newcommand{\fs}{\mathfrak{s}}
\newcommand{\iwa}{\mathfrak{i}}
\newcommand{\fpara}{\mathfrak{p}}
\newcommand{\fp}{\mathfrak{p}}
\newcommand{\fq}{\mathfrak{q}}
\newcommand{\height}{\mathrm{ht}}
\newcommand{\pow}{\mathfrak{o}}

\newcommand{\laur}{F}
\newcommand{\om}{\omega}

\newcommand{\Spec}{\mathrm{Spec}}

\newcommand{\gc}{\nabla}
\newcommand{\fc}{\widehat{\gc}}

\newcommand{\fxnfield}{\C(z)}
\newcommand{\localforms}{\Omega^1_{F/\C}}
\newcommand{\meroforms}{\Omega^1_{\fxnfield/\C}}

\newcommand{\prt}{\mathrm{Part}}
\newcommand{\prn}[2]{\lambda^{{#2},{#1}}}
\newcommand{\prp}[2]{\orb_{#2}^{#1}}
\newcommand{\orbp}[1]{\Orb_{#1}}

\newcommand{\Orb}{\mathrm{Orb}}
\newcommand{\slope}{\mathrm{slope}}

\DeclareMathOperator{\Res}{\mathrm{Res}}
\DeclareMathOperator{\Tr}{Tr}
\DeclareMathOperator{\Lie}{\mathrm{Lie}}

\DeclareMathOperator{\spa}{span}
\newcommand{\Gm}{\mathbb{G}_m}

\newcommand{\hV}{\hat V}

\newcommand{\dzz}{\frac{dz}{z}}
\newcommand{\tdzz}{\tfrac{dz}{z}}

\newcommand{\sdfrac}[2]{\mbox{\small$\displaystyle\frac{#1}{#2}$}}

\newcommand{\ltfrac}[2]{\mbox{\large$\frac{#1}{#2}$}}

\title{Meromorphic connections on the projective line with specified
  local behavior}

\author[D. S. Sage]{Daniel S. Sage}
\address{Department of Mathematics, Louisiana State University, Baton Rouge, LA.}
\email{sage@math.lsu.edu}


\begin{document}
  
\thanks{The author received support from Simons Collaboration Grant
  637367.}
 
\begin{abstract}

  A meromorphic connection on the complex projective line induces
  formal connections at each singular point, and these formal
  connections constitute the local behavior at the singularities.  In
  this primarily expository paper, we discuss the extent to which
  specified local behavior at singular points determines the global
  connection.  In particular, given a finite set of points and a
  collection of ``formal types'' at these points, does there exist a
  moduli space of meromorphic connections with this local behavior,
  and if so, when is this moduli space nonempty or a singleton?  In
  this paper, we discuss variants of these problems (for example, the
  Deligne--Simpson and rigidity problems) as the allowed singularities
  get progressively more complicated: first connections with only
  regular singularities, next connections with additional unramified irregular
  singularities allowed, and finally the general case.

 \end{abstract}

\maketitle

\tableofcontents

\section{Introduction}

Consider a system of first-order linear ordinary differential
equations on the Riemann sphere with rational coefficients:
\begin{equation}\label{ODE} X'+A(z)X=0,
\end{equation}
 where $A(z)$ is an $n\times n$ matrix whose entries are
rational functions.  If the matrix $A(z)$ is not constant, then this
differential equation will have singularities at the finite collection
of points where these matrix entries have poles.  By
expanding the rational functions in Laurent series at each singular
point, one obtains a collection of ``local'' systems of differential
equations which encode the local behavior of the original system.  A
natural question then arises: To what extent does this local behavior
determine the global differential equation?

In order to make this question more precise, it is convenient to give
a geometric reformulation in terms of the equivalent language of
\emph{meromorphic connections}.  Let $V$ be a rank $n$ vector bundle
over the Riemann sphere $\pp$ endowed with a meromorphic connection
$\gc$, i.e., a $\C$-linear derivation $V\to V\otimes\meroforms$ taking
sections of the vector bundle to sections of the corresponding bundle
of one-forms. We will always assume that $V$ is trivializable.  Upon
fixing a trivialization, the connection may be expressed as
$\gc=d+A(z)dz$, where $A(z)\in\fgl_n(\C(z))$ is a matrix of rational
functions.  This is the connection corresponding to our original
differential equation.

Given $p\in\pp$, the local behavior of $\gc$ is governed by the
corresponding \emph{formal connection} $\fc_p$. Let $F=\C(\!(z)\!)$
denote the formal Laurent series, so that $\Delta^\times=\Spec(F)$ is
the formal punctured disk.  A vector bundle over $\Delta^\times$ is
just an $F$-vector space $\hV$, and a formal connection $\fc$ is a
derivation $\hV\to \hV\otimes\localforms$.  One gets the formal
connection $\fc_p$ by expanding the entries of the matrix $A(z)$ as
Laurent series in terms of a local parameter at $p$: $z-p$ if $p\in\C$
and $z^{-1}$ if $p=\infty$.  

As usual, we are only interested in connections (and differential
equations) up to isomorphism, which is given here by \emph{gauge
  equivalence}.  For ODE's, this means that we do not distinguish
between equations which are related by a change of coordinates $Y=gX$,
where $g$ is an invertible $n\times n$ matrix whose coefficients are
regular functions on the base space, namely $\pp$ or $\Delta^\times$.  Accordingly, $g$ is
in $\GL_n(\C)$ or $\GL_n(F)$, depending on whether we are in the
global or local case.  It is an easy calculation to show that our
original ODE \eqref{ODE} becomes\begin{equation}
  Y'+(gA(z)g^{-1}-g'g^{-1})Y=0
\end{equation}
in  the new variables.  For meromorphic ODE's, the nonlinear gauge
term disappears, so that changing coordinates just amounts to
conjugating the original matrix $A(z)$.

In terms of connections, we first observe that the set of trivializations of
a rank $n$ vector bundle over $\Delta^\times$ is a $\GL_n(F)$-torsor.
If $\fc=d+A\,dz$ with respect to the trivialization $\phi$, then the
connection operator for the new trivialization $g\phi$ is obtained by
conjugating the original operator:
\begin{equation} g(d+A\,dz)g^{-1}=d+(gA(z)g^{-1}-g'g^{-1})\,dz.
\end{equation}  
Similar considerations apply to meromorphic connections.

As we have seen, a meromorphic connection $\gc$ on $\pp$ gives rise to
a finite set of singular points and associated formal connections at
these points.  We view the isomorphism classes of these formal
connections as the local behavior at the singularities.  It is now
natural to consider the fiber of this map from global connections to
collections of local connections.  In particular, we are interested in
the following basic questions.

\begin{Three}\label{three}  Fix points  $a_1,\dots,a_k\in\pp$ and formal
  connections $\fc_1,\dots,\fc_k$.
\begin{enumerate}\item{(Existence)} Does there exists a
  meromorphic connection $\gc$, regular on
  $\pp\setminus\{a_1,\dots,a_k\}$, such that for all $i$, the associated
  formal connection $\fc_{a_i}$ is isomorphic to $\fc_i$?
\item{(Uniqueness/Rigidity)} If such a meromorphic connection exists,
  is it unique?  If so, it is called \emph{physically rigid}.
\item{(Moduli spaces)} Is there a moduli space of meromorphic connections with this
  specified local behavior?
\end{enumerate}
\end{Three}

In this paper, we will discuss variants of these problems as the
allowed singularities get progressively more complicated.  A formal
connection is called \emph{regular singular} if there exists a
trivialization for which the connection matrix has a simple pole.
Otherwise, the connection is called \emph{irregular singular}.  

In \S\ref{s:regsing}, we consider meromorphic connections which
only have regular singularities.  In particular, we describe
Crawley--Boevey's solution of the Deligne--Simpson problem for Fuchsian
connections.  In the following section, we turn to a discussion of
connections which also have ``unramified'' irregular singularities.
We explain the construction of moduli spaces of ``unramified framable
connections'' due to Boalch, Hiroe, and Yamakawa~\cite{Boa, HirYam}
and Hiroe's solution of the unramified Deligne--Simpson problem.  In
\S\ref{s:ramsing}, we discuss ramified irregular singular points
with an emphasis on ``toral singularities''.  We construct the moduli
space of framable connections where all ramified singular points are
toral and formulate the \emph{toral Deligne--Simpson problem}.  We then
restrict to a class of meromorphic connection on $\Gm$ called
\emph{Coxeter connections}--certain generalizations of the well-known
Frenkel--Gross and Airy connections.  We describe the solution of the
Deligne--Simpson problem for Coxeter connections; under additional
hypotheses, we determine when they are rigid.  In the final section,
we discuss generalization of these problems to meromorphic
$G$-connections for complex reductive groups $G$.

It is a great pleasure to thank Neal Livesay for many helpful comments
and suggestions.

\section{Connections with only regular singularities}\label{s:regsing}

\subsection{Monodromy of regular singular connections}
It is well-known that the isomorphism class of a regular singular
formal connection is determined by its \emph{monodromy}.  The
monodromy is easy to compute as long as one puts a very mild
restriction on the leading term of the connection.  Consider the
formal connection $\fc=d+(B_0+B_1 z+\dots)\dzz$, where
$B_k\in\fgl_n(\C)$.  It is called \emph{nonresonant} if no two
eigenvalues of $B_0$ differ by a nonzero integer.

\begin{lem}\label{regsinggauge}  If the connection $\fc$ is nonresonant, then it is
  formally gauge equivalent to $d+B_0\dzz$.
\end{lem}

We briefly recall the proof.  It suffices to  find $g=g_0+g_1 z+g_2
z^2+\dots\in\GL_n(\C[\![z]\!])$ with $g_i\in\fgl_n(\C)$ such that $B_0/z=g(\sum_{k=0}^\infty
B_kz^{k-1})g^{-1}-g'g^{-1}$.  The power series $g$ is invertible as
long as $g_0\in\GL_n(\C)$.   Multiplying through by $zg$ on the right
and rearranging gives the following differential equation
for $g$:
\begin{equation} zg'=g(\sum_{k=0}^\infty
B_kz^k)-B_0g.
\end{equation}   If we set $g_0=I$, the $z^0$ coefficient on
both sides is $0$.  Equating the coefficients of $z^k$ gives
\begin{equation}\label{regsingrec} (B_0+kI)g_k-g_kB_0=\sum_{i=0}^{k-1}g_i B_{k-i}.
\end{equation}
We now consider the linear map $f_{C,D}:\fgl_n(\C)\to\fgl_n(\C)$ given
by $X\mapsto CX-X D$, where $C,D\in\fgl_n(\C)$.  It is a standard fact
from matrix theory that this map is an isomorphism if and only if $C$
and $D$ do not have a common eigenvalue~\cite[Theorem 4.1]{Was}.
Applying this result to $C=B_0+kI$, $D=B_0$, and using the fact that
$B_0$ is nonresonant, we see that we can solve the equations
\eqref{regsingrec} recursively to obtain the desired gauge change.

In the nonresonant case, we are now reduced to computing the monodromy
of $d+B_0\dzz$ or equivalently, of the ODE $X'+z^{-1}B_0 X=0$.  A
fundamental matrix of this equation is given by a branch of
$z^{-B_0}=e^{-B_0\log{z}}$.  Analytically continuing a local solution
at $1$ along a counterclockwise loop around $0$ gives a new
fundamental solution which differs from the first by the monodromy
matrix $e^{-2\pi i B_0}$.  The monodromy is only determined up to
conjugation by a constant matrix; indeed, a constant gauge change of
$d+B_0\dzz$ acts by conjugation on $B_0$ and hence on the monodromy.
Thus, the isomorphism class is given by the conjugacy class of
$e^{-2\pi i B_0}$ in $\GL_n(\C)$.  In particular, the adjoint orbit of
the leading term determines the isomorphism class of the formal
connection, though different orbits can lead to isomorphic formal
connections.  We will view this adjoint orbit (or any of its
representatives) as a \emph{formal type} of $\fc$.

\subsection{Fuchsian connections}  

It is difficult to study our three questions for arbitrary meromorphic
connections with only regular singularities.  Instead, we will
restrict attention to \emph{Fuchsian connections}.  These are
meromorphic connections which have a global representation in which
all singularities are manifestly simple poles.  We will assume without
loss of generality that $\infty$ is a regular point.

\begin{definition}  A meromorphic connection $\gc$ is called \emph{Fuchsian} with singular
  points $a_1,\dots,a_k$ if there exists a trivialization in which
\begin{equation} \label{fuchs}\gc=d+\sum_{i=1}^k\frac{A_i}{z-a_i}dz,
\end{equation}
where $A_i\in\fgl_n(\C)$.
\end{definition}
Since we are assuming that $\infty$ is not a singular point, the
residue theorem implies that $\sum_{i=1}^k A_i=0$.

The local behavior of a Fuchsian connection is determined by a
collection of nonresonant adjoint orbits $\orb_1,\dots,\orb_k$, where
$\orb_i$ is the formal type at $a_i$. A Fuchsian connection has
local behavior determined by these formal types if $A_i\in\orb_i$ for
all $i$; we let $\sM(\orb_1,\dots\orb_k)$ denote the corresponding
moduli space.  (The moduli space does not depend on the specific $k$
points in $\pp$, so we omit the points from the notation.)

\begin{prop} \label{p:fuchsmodspace} Given nonresonant adjoint
orbits $\orb_1,\dots,\orb_k\subset\fgl_n(\C)$, the corresponding moduli
space of Fuchsian connections is given by
\begin{equation}\label{e:fuchsmodspace}\sM(\orb_1,\dots\orb_k)=\{(A_1,\dots,A_k)\in\orb_1\times\dots\times\orb_k\mid\sum_{i=1}^k A_i=0\}/\GL_n(\C).
\end{equation}
\end{prop}

We remark that the form of the moduli space is very natural from the
point of view of symplectic geometry.  Adjoint orbits can be viewed as
symplectic manifolds via the identification with coadjoint orbits
given by the trace form.  The action of $\GL_n(\C)$ on the symplectic
manifold $\orb_1\times\dots\times\orb_k$ is Hamiltonian with moment map
$\mu:(A_1,\dots,A_k)\mapsto \sum_{i=1}^k A_i$.  The moduli space is
thus the Hamiltonian reduction 
\begin{equation} (\orb_1\times\dots\times\orb_k)\sslash_0\GL_n(\C):=\mu^{-1}(0)/\GL_n(\C).
\end{equation}

\subsection{The Deligne--Simpson problem}

While it is easy to find the explicit form for the moduli
space~\eqref{e:fuchsmodspace}, it is more challenging to understand
its geometry.  In fact, it is quite difficult to determine when the
moduli space is nonempty except for very small values of $k$ (or
$n=1$).  If there is only one singular point, then the moduli space is
a singleton (consisting of the trivial connection $d$) if
$\orb_1=\{0\}$ and is empty otherwise.  If $k=2$, then the moduli
space is empty unless $\orb_2=-\orb_1$, in which case it is again a
singleton.

To go further, we restrict attention to the \emph{stable moduli
  space} consisting of irreducible Fuchsian connections.  Explicitly,
this means that the residue matrices do not have any nonzero proper
invariant subspaces in common.  Note that for $k=2$ (and $n\ge 2$),
the stable moduli space is always empty.  We can now state the
(additive) Deligne--Simpson problem:

\begin{quote} \emph{Given adjoint orbits $\orb_1,\dots,\orb_k$,
    determine whether there exists an irreducible $k$-tuple
    $(A_1,\dots,A_k)$ with $A_i\in\orb_i$ satisfying $\sum_{i=0}^k
    A_i=0$.  Equivalently, determine whether there exists an
    irreducible Fuchsian connection with residues in the prescribed orbits~\cite{Kostov03}.}
\end{quote}

The original Deligne--Simpson problem was a multiplicative version
involving the equation $M_1\dots M_k=1$ with the matrices $M_i$ lying
in fixed conjugacy classes of $\GL_n(\C)$~\cite{Simpson}.  In terms of
connections, this is the existence problem for irreducible Fuchsian
connections with given monodromy classes, and a solution for the
additive problem gives rise to a solution for the multiplicative one.
The additive Deligne--Simpson problem was formulated by Kostov, who
solved it under various genericity
hypotheses~\cite{Kostov02,Kostov03,Kostov04}.  A complete solution was
given by Crawley-Boevey by reinterpreting the problem in terms of
representations of quivers~\cite{CB03}.

We now sketch Crawley-Boevey's solution. 

\noindent {\bf Step 1. }
One begins by associating three objects to the collections of orbits:
a quiver $Q$, say with vertex set $I$, a dimension vector
$\alpha\in\Z_{\ge 0}^I$, and a deformation vector $\lambda\in\C^I$.
The quiver only depends on the degrees of the minimal polynomials of
the orbits while the definition of the dimension vector also involves
the sizes of the Jordan blocks.  Only the deformation vector depends
on the specific eigenvalues appearing in the
orbits.  

Let $p_i(x)=\prod_{j=1}^{d_i} (x-\eta_{ij})$ be the minimal polynomial
of $\orb_i$.  The quiver $Q$ has a sink labeled by $0$ and a path of
length $d_i-1$ leading into it for each orbit.  The vertices on the
$i$th path are labeled by $[i,j]$ for $1\le j< d_i$.\footnote{If
  $\orb_i$ is a scalar matrix, then the path has length $0$, and the
  orbit does not contribute any additional vertices to the quiver.}

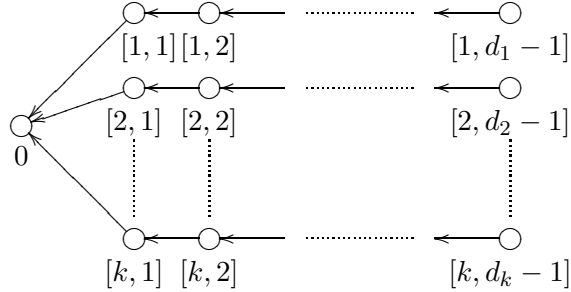
 \[
            \begin{xy}
            (-10,0) *++!U{0} *\cir<4pt>{}="A",
            (5,15) *++!U!L(0.3){[1,1]} *\cir<4pt>{}="B",
            (15,15) *++!U{[1,2]} *\cir<4pt>{}="C",
            (55,15) *++!U{[1,d_{1}-1]} *\cir<4pt>{}="D",
            (5,5) *++!U{[2,1]} *\cir<4pt>{}="E",
            (15,5) *++!U{[2,2]} *\cir<4pt>{}="F",
            (55,5) *++!U{[2,d_{2}-1]} *\cir<4pt>{}="G",
            (5,-15) *++!U{[k,1]} *\cir<4pt>{}="H",
            (15,-15) *++!U{[k,2]} *\cir<4pt>{}="I",
            (55,-15) *++!U{[k,d_{k}-1]} *\cir<4pt>{}="J",
            \ar@{<-} "A";"B"
            \ar@{<-} "A";"E"
            \ar@{<-} "A";"H"
            \ar@{<-} "B";"C"
            \ar@{<-} "C";(25,15)
            \ar@{.} (28,15);(42,15)
            \ar@{<-} (45,15);"D"
            \ar@{<-} "E";"F"
            \ar@{<-} "F";(25,5)
            \ar@{.} (28,5);(42,5)
            \ar@{<-} (45,5);"G"
            \ar@{<-} "H";"I"
            \ar@{<-} "I";(25,-15)
            \ar@{.} (28,-15);(42,-15)
            \ar@{<-} (45,-15);"J"
            \ar@{.} (5,-2);(5,-12)
            \ar@{.} (15,-2);(15,-12)
            \ar@{.} (55,-2);(55,-12)
        \end{xy}
        \]

        To define the dimension vector $\alpha$, we take $C_i\in\orb_i$,
        and let $r_{ij}$ be the rank of
        $\prod_{\ell=1}^{j} (C_i-\eta_{i\ell})$ for $0\le j\le d_i$.
        We then set $\alpha_0=n$ and $\alpha_{[i,j]}=r_{i,j-1}-r_{ij}$.
        Finally, the deformation vector is defined by
        $\lambda_0=-\sum_{i=1}^k \eta_{i1}$ and
        $\lambda_{[i,j]}=\eta_{ij}-\eta_{i,j+1}$.

\noindent {\bf Step 2. }    Let  $\overline{Q}$ be the doubled quiver of
$Q$.  This is a quiver on the same vertex set $I$, but for each arrow
$a$ from $i$ to $i'$ in $Q$, there is a new arrow $a^*$ from $i'$ to
$i$.  One can now consider the representations $(V_i)_{i\in I}$ of $\overline{Q}$
with dimension vector $\alpha$ (i.e., $\dim V_i=\alpha_i$) and which satisfy the relations
\begin{equation}\label{defcond}
  \sum_{h(a)=i}aa^*-\sum_{t(a)=i}a^*a=\lambda_i 1_{V_i},
\end{equation}
where $h(a)$ and $t(a)$ denote the head and tail of the arrow
$a$. Equivalently, these are the representations of dimension $\alpha$
of the deformed preprojective algebra consisting of the path algebra
$\C \overline{Q}$ modulo the ideal generated by the given relations.

Crawley-Boevey shows that there is a bijection between
irreducible solutions of the the equation $\sum_{i=1}^k A_i=0$ with
$A_i\in\orb_i$ and irreducible representations of $\overline{Q}$ with
dimension vector $\alpha$ and satisfying the deformation
conditions~\eqref{defcond}.  In fact, this bijection induces a
symplectic isomorphism between the corresponding moduli spaces.

\noindent {\bf Step 3. }  It remains to determine when such
representations of $\overline{Q}$ exist.  The solution to this problem
is given in terms of the Kac--Moody Lie algebra $\fg_Q$ associated to
the quiver $Q$.  A symmetric matrix with integer coefficients $C$ is
called a \emph{symmetric generalized Cartan matrix} if its diagonal
entries are all $2$ and its off-diagonal entries are nonpositive.
Such a matrix gives rise to a corresponding Kac--Moody Lie algebra.  In
particular, one gets $\fg_Q$ from the matrix $C_Q$ defined
by \begin{equation} (C_Q)_{ij}=2\delta_{ij}-\text{ $\#\{$edges from
    $i$ to $j$\}}.
\end{equation}
The root system $R$ of $\fg_Q$ is a certain subset of $\Z^I$.  A root
is called positive if it lies in $\Z_{\ge 0}^I$, and $R=R_+\sqcup
-R_+$, where $R_+$ is the set of positive roots.

Let $R_+^\lambda=\{\beta\in R_+\mid \beta\cdot\lambda=0\}$.  Let
$p:\Z^I\to \Z$ be the function defined by $p(\beta)=1-(1/2)\beta^t
C_Q\beta$.  If $\beta$ is a root, then $p(\beta)\ge 0$.  We now let
$\Sigma^\lambda$ be the set of $\beta\in R_+^\lambda$ such that if
$\beta=\sum_{j=1}^\ell\gamma_j$ with$\gamma_j\in R_+^\lambda$ and 
$\ell\ge 2$, then $p(\beta)>\sum _{j=1}^\ell p(\gamma_j)$.

Crawley-Boevey showed that there exists an irreducible representation
of $\overline{Q}$ of dimension $\alpha$ satisfying~\eqref{defcond} if
and only if $\alpha\in \Sigma^\lambda$.  In particular:

\begin{thm}[\cite{CB03}]  There is an irreducible Fuchsian
  connection with residues in the $\orb_i$'s if and only if $\alpha\in
  \Sigma^\lambda$.
\end{thm}

\subsection{Rigidity} There are several different notions of rigidity
for irreducible Fuchsian connections.  One can ask whether such a
connection is physically rigid, i.e., uniquely determined by the
formal isomorphism classes at its singular points.  One can also ask
whether the moduli space of Fuchsian connections with the given formal
types is a singleton.  In this section, we will call an irreducible
Fuchsian connection rigid if the corresponding stable moduli space is
a singleton.  Note that a rigid Fuchsian connection may not be
physically rigid; for example, there could be a reducible Fuchsian
connection with the same formal types.

In~\cite{CB03}, Crawley-Boevey also determined when an irreducible
Fuchsian connection is rigid.  To explain this, we need to recall a
bit more about the root system $R$ of $\fg_Q$.

The coordinate vectors $e_i\in\Z^I$ lie in $R$; they are called the
simple roots of $R$.  The symmetric matrix $C_Q$ defines a bilinear
form on $\Z^I$, and so one can define the simple reflection
$s_i:\Z^I\to \Z^I$ in the root $e_i$ given by
$s_i(\beta)=\beta-(\beta^t C_Q e_i)e_i$.  The Weyl group $W$ of
$\fg_Q$ is the group generated by the $s_i$'s.  A root $\beta$ is
called \emph{real} if $\beta=w e_i$ for some $w\in W$ and $i\in I$;
otherwise, it is called \emph{imaginary}.  For example, if the quiver
$Q$ is a Dynkin quiver (i.e., the underlying undirected graph is a
simply-laced Dynkin diagram), then all roots are real.  The real roots
can also be characterized as the roots satisfying $p(\beta)=0$.

It turns out that if there exists an irreducible Fuchsian connection
for the adjoint orbits $\orb_1,\dots,\orb_k$, then it is rigid if and
only if $\alpha$ is a real root.  In fact, one can say more.

\begin{thm}[\cite{CB03}]  Suppose that  $\alpha\in
  \Sigma^\lambda$.  Then, the stable moduli space
  $\sM^{\mathrm{st}}(\orb_1,\dots,\orb_k)$ is a singleton if $\alpha$
  is real and is infinite if $\alpha$ is imaginary.
\end{thm}
 If $\alpha$ is real, then it is more straightforward to check whether
 $\alpha\in\Sigma^\lambda$: if $\alpha$ is a real positive root, then $\alpha\in\Sigma^\lambda$ if and only if $\alpha\cdot\lambda=0$ and
for any nontrivial decomposition of $\alpha$ into a sum of positive
roots $\alpha=\sum_{j=1}^\ell \beta_j$, $\beta_k\cdot\lambda\ne 0$ for
some $k$.

\subsection{Examples}

\begin{example} Suppose that $k=2$ (with $n\ge 2$ arbitrary), so that
  there are two singular points.  In this case, the underlying graph
  of the quiver is the Dynkin quiver $A_{d_1+d_2-1}$.  Every positive
  root of this quiver is a sum of consecutive simple roots, all with
  multiplicity $1$.  However, $e_0$ appears in $\alpha$ with
  multiplicity $n\ge 2$, so we recover the elementary fact that the
  stable moduli space is empty in this case.
\end{example}

\begin{example} We now consider the case of three singular points on a
  rank $2$ connection, i.e., $k=3$, $n=2$~\cite{CB03}.  We will show
  the following:
  \begin{equation*}
  \begin{aligned}
    \text{An irreducible}& \text{ solution exists}\\
    \text{(and is}& \text{ unique)}
  \end{aligned}
  \quad \iff\quad
  \begin{aligned} &\text{(a)    $\,\,\sum_{i=1}^3\Tr(\orb_i)=0$;}\\
    &\text{(b) $\,$   the sum of $3$ eigenvalues, one from each $\orb_i$, is nonzero.}
  \end{aligned}
\end{equation*}

We denote the eigenvalues of the (nonscalar) adjoint orbits $\orb_i$
by $\eta_{11},\eta_{12}$.  Since the minimal polynomial of each orbit
has degree $2$, the corresponding quiver has its underlying graph
given by the $D_4$ Dynkin diagram.  The dimension vector is
$\alpha=(2,1,1,1)$, so the quiver representations are given by the
diagram below.  The deformation vector is
           \begin{equation*}\lambda=(-(\eta_{11}+\eta_{21}+\eta_{31}), \eta_{11}
           -\eta_{12}, \eta_{21} -\eta_{22}, \eta_{31} -\eta_{32}).
         \end{equation*}
         
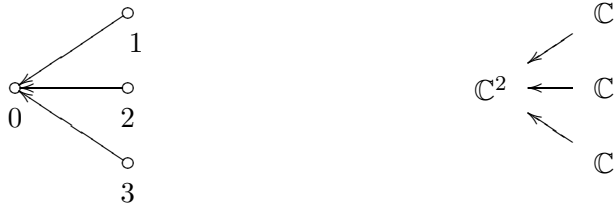
\[
            \begin{xy}
            (-10,0) *++!U{0} *\cir<2pt>{}="A",
            (5,10) *++!U!L(0.3){1} *\cir<2pt>{}="B",
             (5,0) *++!U{2} *\cir<2pt>{}="E",
             (5,-10) *++!U{3} *\cir<2pt>{}="H",
            \ar@{<-} "A";"B"
            \ar@{<-} "A";"E"
            \ar@{<-} "A";"H"
          \end{xy}\qquad\qquad\qquad\qquad\qquad
          \begin{xy}
            (170,0) *++{\mathbb{C}^2}="A",
            (185,10) *++{\mathbb{C}}="B",
             (185,0) *++{\mathbb{C}}="E",
            (185,-10) *++{\mathbb{C}}="H",
            \ar@{<-} "A";"B"
            \ar@{<-} "A";"E"
            \ar@{<-} "A";"H"
          \end{xy}
           \]

           The $D_4$ root system has $12$ positive roots: $e_0$,
           $e_1$, $e_2$, $e_3$, $e_0+e_1$, $e_0+e_2$,$e_0+e_3$,
           $e_0+e_1+e_2$, $e_0+e_1+e_3$, $e_0+e_2+e_3$,
           $e_0+e_1+e_2+e_3$, and $2e_0+e_1+e_2+e_3$.  Thus, $\alpha$
           is a (real) root.  The condition that
           $\alpha\cdot\lambda=0$ is just the trivial necessary
           condition that if there exist matrices $A_i\in\orb_i$
           summing to $0$, then their traces sum to $0$.  Given
           $J\subset\{1,2,3\}$, let $\beta_J=e_0+\sum_{i\in J} e_i$.
           Note that $\beta_J\cdot\lambda=-(\sum_{i\in J}
           \eta_{i2}+\sum_{i\notin J} \eta_{i1})$, so the stated
           condition on the sum of $3$ eigenvalues means precisely
           that $\beta_J\cdot\lambda\ne 0$ for all $J$.

Every decomposition of $\alpha$ into the sum of two positive roots is
of the form $\alpha=\beta_J+\beta_{J^c}$, so for an irreducible
solution to exist, we need either $\beta_J\cdot\lambda\neq 0$ or
$\beta_{J^c}\cdot\lambda\neq 0$.  However, under the assumption that
$\alpha\cdot\lambda=0$, these are equivalent.  Thus,
$\alpha\in\Sigma^{\lambda}$ implies that $\beta_J\cdot\lambda\neq 0$
for all $J$.  On the other hand, if this condition on the $\beta_J$'s
is satisfied, then $\alpha\in\Sigma^{\lambda}$.  Indeed, any
nontrivial decomposition of $\alpha$ into a sum of positive roots
involves at least one $\beta_J$.

      \end{example}

\section{Connections with unramified irregular singular points}\label{s:unramsing}

\subsection{The Levelt--Turrittin normal form}

We now shift our attention from connections with only regular
singularities to connections where irregular singular points are allowed.
First, we need to understand isomorphism classes of irregular singular
formal connections and to get an appropriate notion of formal type.
To do this, we recall that there is an analogue of Jordan canonical
form for formal connections called the \emph{Levelt--Turrittin (LT)
  normal form}.  This allows one to ``upper-triangularize'' a formal
connection at the cost of possibly introducing roots of the formal
parameter $z$.

\begin{thm}[Levelt--Turrittin normal form \cite{Tur55,Lev75}]  Any
  formal connection is isomorphic to a connection of the form
  \begin{equation}\label{LT} d+(D_rz^{-r/b}+D_{r-1}z^{-(r-1)/b}+\dots+ D_1z^{-1}+D_0+N)\sdfrac{dz}{z},
\end{equation} where the $D_i$'s are diagonal, $D_r\ne 0$, $N$ is
strictly upper triangular and commutes with each $D_i$, and $b\in\N$.
\end{thm} 

Note that the LT normal form does not contain any terms in degrees
above the residue term, and it is diagonal except possibly for a
nilpotent component to the residue.  The rational number $r/b$ is an
invariant of the formal connection called its \emph{slope}.  It turns
out that $\fc$ is irregular singular if and only if its slope is
positive.  Moreover, the denominator of the slope in lowest form is at
most $n$. (This does not imply that one can always take $b\le n$.)

There is a more geometric formulation of the Levelt--Turritin form as
the \emph{slope decomposition} of a formal connection.  Let
$\Delta_b^\times=\Spec(\C(\!(z^{1/b})\!))$; it is a degree $b$
ramified cover of the formal punctured disk $\Delta^\times$.  It
follows from this theorem that after passing to this ramified cover,
the pullback connection is isomorphic to a finite direct sum
\begin{equation} \bigoplus(L_i\otimes M_i,\fc_{L_i}\otimes\fc_{M_i}),
\end{equation}
where $L_i$ is a line bundle and $\fc_{M_i}$ is regular singular.  One
can now consider the collection of nonnegative rational numbers
$\slope(\fc_{L_i})$, each taken with multiplicity $\dim M_i$; these
are called the \emph{slopes} of $\fc$.  More concretely, they are the
negatives of the valuations of the eigenvalues of the matrix of the LT
form.\footnote{The slopes corresponding to $0$ eigenvalues are taken
  to be $0$.}  These are also invariants of $\fc$, and the slope of
$\fc$ is the maximum of the slopes, i.e.,
$\slope(\fc) = \max_i (\slope(\fc_{L_i}))$.  The \emph{irregularity} of
$\fc$ is the sum of the slopes; it is a nonnegative integer which is
$0$ if and only if $\fc$ is regular singular~\cite{Katz87}.

\begin{definition} \begin{enumerate} \item A formal connection $\fc$ is called
  \emph{unramified} if the LT form (or the slope decomposition) may be
  obtained without passing to a ramified cover (i.e. if one may take
  $b=1$).  In particular, all of the
    slopes are integers.
  \item The set of  \emph{unramified formal types of slope $r$} is the
    collection of matrix one-forms $(D_{-r}z^{-r}+\dots+ D_0+N)\dzz$ with
the $D_i$'s diagonal, $D_{-r}\ne 0$, and $N$ strictly upper triangular and commuting
with each $D_i$.  If $r=0$, we further assume that $D_0$ is nonresonant.
\end{enumerate}
\end{definition}

A regular singular formal connection is an unramified formal
connection of slope (or irregularity) $0$.  The notion of a regular
singular unramified formal type is slightly different, but equivalent
to the regular singular formal types used in the previous section.  We
will continue to use the orbit notation $\orb$ for regular singular
formal types.

It will be helpful in our discussion of the unramified Deligne--Simpson
problem below to express unramified formal types in terms of the
simultaneous eigenspaces of the purely irregular part.  If $\ftype$ is
an unramified formal type of rank $n$, then there exists a
decomposition $V=V_1\oplus\dots\oplus V_\ell$ such that
\begin{equation}\label{ltresform} \ftype=\bigg(\sum_{j=1}^\ell \left(q_j(z)
    I_{V_j}+R_j\right)\bigg)\sdfrac{dz}{z},
\end{equation}
where $q_j(z)\in z^{-1}\C[z^{-1}]$ with $q_j(z)\ne q_{j'}(z)$ for
$j\ne j'$ and 
$R_j\in\fgl(V_j)$.  Note that if $\ftype$ is regular singular, then
$\ell=1$ and $q_1=0$.

\subsection{Unramified framable connections}

As in the regular singular case, we will not consider all meromorphic
connections with specified unramified behavior, but rather will limit
ourselves to an irregular version of Fuchsian connections called
\emph{framable connections}.

Fix singular points $a_1,\dots,a_k$ and
unramified formal types $\ftype_1,\dots,\ftype_k$ with slopes
$r_1,\dots,r_k$.  For convenience, we will assume that $\infty$ is a
regular point.    We will restrict attention to connections of the form 
\begin{equation}\label{eq:pp} d+\bigg(\sum_{i=1}^k \sum_{\nu=0}^{r_i}
  \frac{A^{(i)}_{\nu}}{(z-a_i)^{\nu}}\bigg)\sdfrac{dz}{z}
\end{equation}
for which $A^{(i)}_{r_i}$ is nonnilpotent if $r_i>0$.  We call these
connections \emph{framable}.  Since $\infty$
is regular, the residue theorem implies that
$\sum_{i=1}^k A^{(i)}_{0}=0$. Such connections clearly have a singular
point of slope $r_i$ at $a_i$ and have no other singular points.

We need to clarify what it means for a framable connection to have
unramified formal type $\ftype_i$ at $a_i$.  Recall that for Fuchsian
connections, the principal part reduces to the residue term, and we
require the residue to be in the specified adjoint orbit.  The
analogue of this for irregular unramified formal types is that the
principal part of $\gc$ at $a_i$ must lie in the orbit of $\ftype_i$
under a certain action of the group $\GL_n([\![z-a_i]\!])$.

Let $\cB_r=\{(B_{-r} z^{-r}+\dots + B_0)\dzz\mid B_i\in\fgl_n(\C)\}$
denote the space of principal parts of order at most $r$.  If we let
$\pow=\C[\![z]\!]$ be the ring of formal power series, then the group
$\GL_n(\pow)$ acts on $\cB_r$ by conjugation followed by
truncation at the residue term.  Note that this action factors through
the finite-dimensional group $\GL_n(\pow/z^{r+1} \pow)$.

One can interpret truncated conjugation as a coadjoint action.  Let
$\fgl_n(F)^\vee$ be the space of smooth $\C$-linear functionals on
$\fgl_n(F)$.  Here, smooth means that the functional vanishes on a
nonempty bounded open subgroup.  Any one-form $B\dzz$ with
$B\in\fgl_n(F)$ defines a smooth functional on $\fgl_n(F)$ via
$A\mapsto\Res\Tr(AB\dzz)$.  In fact, if we identify $B\dzz$ with the
corresponding functional, the map $B\to B\dzz$ is an isomorphism
$\fgl_n(F)\to \fgl_n(F)^\vee$.  It is easy to see that
$\fgl_n(\pow)^\perp=z \fgl_n(\pow)\dzz$, so
$\fgl_n(\pow)^\vee=\fgl_n(\C[z^{-1}])\dzz$.  Thus, $\GL_n(\pow)$ acts
on $\fgl_n(\C[z^{-1}])\dzz$ via the coadjoint action with $\cB_r$ as a
subrepresentation, and this action corresponds to truncated
conjugation.

If $\ftype$ is an unramified formal type of slope $r$, we call the
orbit $\orb_{\ftype}\subset\cB_r$ under this action the
\emph{truncated orbit} of $\ftype$.  Since it is a coadjoint orbit, it
has a natural symplectic structure.  Moreover, the $\GL_n(\C)$-action
on $\orb_{\ftype}$ is Hamiltonian with moment map given by the
residue.   If $\ftype$ has slope $0$,
$\orb_{\ftype}$ may be identified with the usual adjoint orbit of
$\ftype/\dzz$.

The framable connection \eqref{eq:pp} has formal type $\ftype_i$ at
$a_i$ if $\sum_{\nu=0}^{r_i}A^{(i)}_{\nu}z^{-\nu}\dzz\in\orb_{\ftype_i}$.
Accordingly, we get the following expression for the moduli space of
framable connections with the given unramified formal types.
\begin{thm}[\cite{Boa,HirYam}] \label{unrammodspace} The
  moduli space of framable connections with unramified formal types $\ftype_i$ at
  $a_i$ is given by
  \begin{equation}\label{unrammod}\sM(\ftype_1,\dots,\ftype_k)=\{(A_1,\dots,A_k)\in\orb_{\ftype_1}\times\dots\times
    \orb_{\ftype_k}\mid \sum_{i=1}^k \Res A_i=0\}/\GL_n(\C).
  \end{equation}
  Equivalently, it is given by Hamiltonian reduction:
  \begin{equation}(\orb_{\ftype_1}\times\dots\times
    \orb_{\ftype_k})\sslash_0\GL_n(\C).
  \end{equation}
\end{thm}

Note that if all the singular points are regular, we recover the
moduli space of Fuchsian connections with specified formal types given
in Theorem~\ref{p:fuchsmodspace}.  Just as in the Fuchsian case, the
location of the singular points is irrelevant.

\begin{rmk} Boalch constructed this moduli space under the additional
  hypothesis that the leading term of every irregular unramified
  formal type is regular semisimple~\cite{Boa} while Hiroe and
  Yamakawa proved the theorem in general~\cite{HirYam}.
\end{rmk}

We end this section by working out this moduli space in the simplest
non-Fuchsian case.  This is a generalization of an example that
appeared in~\cite{Sa17}.
\begin{example}\label{unrameg}  Consider the case of two singular points on a
  rank $2$ connection, i.e., $k=2$, $n=2$, with formal types $\ftype$ unramified
  of slope $1$ and $\orb$ regular singular.  There are two
  possible forms for $\ftype$, depending on whether or not it is diagonalizable:
  \begin{equation}\ftype=\begin{cases}
      \left(\diag(a,b)z^{-1}+\diag(c,d)\right)\dzz &\text{$a,b$ not both $0$;}\\
      \left((az^{-1}+c)I+e_{12}\right)\dzz & a\ne 0.
    \end{cases}
  \end{equation}
 The moduli space is given by
 \begin{equation}\sM(\ftype,\orb)=\{(\alpha,Y)\in\Ad^*(\GL_2(\pow))(\ftype)\times\orb\mid
   \Res(\alpha)+Y=0\}/\GL_2(\C).
 \end{equation}

 In order to compute the truncated orbit, we will use the fact that $\GL_2(\pow)=\GL_2(\C)\ltimes(1+t\fgl_2(\pow))$.

  \subsubsection*{Case 1: $a\ne b$} 
   We first observe that
  $\Ad^*(1+z\fgl_2(\pow))(\ftype)=\left\{\left(\begin{smallmatrix}az^{-1}+c&u\\v&bz^{-1}+d
    \end{smallmatrix}\right)\dzz\right\}$ with $u,v\in\C$ arbitrary. Indeed,
  if $X,X'\in\fgl_2(\C)$, then $(1+zX)X'(1+zX)^{-1}\equiv
  X'+z[X,X']\pmod{z^2}$, and the claim follows since
  $\ad(\diag(a,b))(\fgl_2(\C))$ consists of the off-diagonal matrices.  Thus,
\begin{equation}\label{boalchorbit}
   \Ad^*(\GL_2(\pow))(\ftype)=\Ad^*(\GL_2(\C))\left\{\begin{pmatrix}a z^{-1}+c&u\\v&b
    z^{-1}+d
   \end{pmatrix}\dzz
   \Big| u,v\in\C\right\}.
\end{equation}

It follows from \eqref{boalchorbit} that every $\GL_2(\C)$-orbit of
pairs $(\alpha, Y)$ has a representative with $\alpha$
  in the standard
  form $\left(\begin{smallmatrix} az^{-1}+c&u\\v&bz^{-1}+d
  \end{smallmatrix}\right)\dzz$ for some $u,v\in\C$.  We view the pair
$(u,v)$ as coordinates on the set of standard $\alpha$'s.  Since the
diagonal torus $T$ is the
stabilizer of the diagonal component, hence preserves the set of
standard forms, it follows that the moduli space is
the set of $T$-orbits of pairs $(\alpha,Y)$ with $\alpha$ in
standard form.  We claim that 
\begin{equation}|\sM(\ftype,\orb)|=\begin{cases} 3,&\text{$c\ne d$, $\det(\orb)=cd$ and
      $\Tr(\orb)=-c-d$;}\\ 2,&\text{$c=d$, $\det(\orb)=cd$ and
      $\Tr(\orb)=-c-d$;}\\ 1,&\text{$\det(\orb)\ne cd$ and
      $\Tr(\orb)=-c-d$;}\\
    1,&\text{$c=d$ and $\orb=\{-cI\}$;}\\
0,&\Tr(\orb)\ne -c-d.
\end{cases}
\end{equation}

The torus $T$ acts on the set of standard $\alpha$'s via
$\diag(s,t)\cdot (u,v)=(st^{-1}u,ts^{-1}v)$.  Thus, the $T$-orbits are
parameterized by the pairs $(u,1)$ with $u\in\C$, $(1,0)$, and
$(0,0)$.  The residue of a standard $\alpha$ can have minimal
polynomial of degree $1$ only if $c=d$, and then one needs
$(u,v)=(0,0)$.  Accordingly, we may now assume that the minimal
polynomial of $\orb$ is the characteristic polynomial.  Since
$\Res(\alpha)+Y=0$, their traces also sum to zero, implying that the
moduli space is empty unless $\Tr(\orb)= -c-d$.  If $\det(\orb)\ne
cd$, then there is a unique $T$-orbit representative which equates
$\det(\Res(\alpha))$ and $\det(\orb)$, namely $(cd-\det(\orb),1)$.
Thus, in this case, the moduli space is a singleton.  On the other
hand, if $\det(\orb)= cd$, the three orbit representatives $(0,1)$,
$(1,0)$, and $(0,0)$ all give rise to a standard $\alpha$ with
$\det(\Res(\alpha))=\det(\orb)$.  If $c\ne d$, these give the $3$
points in the moduli space.  If $c=d$, then only the first two give
rise to a nonscalar adjoint orbit, so there are only $2$ points in the
moduli space.

\subsubsection*{Case 2: $a=b$}  In this case,
 $1+z\fgl_2(\pow)$ fixes $\ftype$, so
 $\Ad^*(\GL_2(\pow))(\ftype)=\Ad^*(\GL_2(\C))(\ftype)$ is a single $\GL_2(\C)$-orbit.
 It follows that $\sM(\ftype,\orb)$ is empty unless
 $\Res(\ftype)\in-\orb$, in which case the moduli space is a singleton.
\end{example}

\subsection{The unramified Deligne--Simpson problem}

As in the Fuchsian case, it is not easy to understand the moduli space
\eqref{unrammod}.  It can be quite badly behaved; in fact, it does not
even have to be Hausdorff~\cite{Boa}.  However, as in the Fuchsian
case, the moduli space is much nicer if one restricts to irreducible
framable connections.  The framable connection \eqref{eq:pp} is
irreducible if the collection of matrices $\{A^{(i)}_{\nu}\}$ have no
common nontrivial proper invariant subspace.  Hiroe has shown that the
stable moduli space consisting of the irreducible framable connections
is a connected complex manifold if it nonempty~\cite{Hiroe}.  The
unramified Deligne--Simpson problem is simply the determination of when
this stable moduli space is nonempty.

\begin{example}  We return to the situation of Example~\ref{unrameg}.
  If $a\ne b$, then $\sM^{\mathrm{st}}(\ftype,\orb)$ is empty except
    when $\det(\orb)\ne cd$ and
      $\Tr(\orb)=-c-d$, in which case it is a singleton.  Indeed, in all the cases where
      $\det(\orb)=cd$, there is a representative $(\alpha,Y)$ of each
      solution with $\alpha=\ftype$ (hence diagonal) and $Y$ either
      upper or lower triangular.  If $\det(\orb)\ne cd$, then the
      unique solution has a representative with $\alpha$ having
      leading term $\diag(a,b)z^{-1}\dzz$ and $Y=-\left(\begin{smallmatrix} c&u\\1&d
  \end{smallmatrix}\right)$ with $u\ne 0$.  This means that the two
eigenspaces of $\diag(a,b)$ are not invariant  under $Y$.   Finally,
if $a=b$, then it is trivial that the stable moduli space is empty as
the $3$ relevant matrices are $aI$, $\Res(\alpha)$ and $Y=-\Res(\alpha)$.
 \end{example}

Hiroe has given a complete solution of the unramified Deligne--Simpson
problem, building on earlier work of Boalch~\cite{Boa} and Hiroe and
Yamakawa~\cite{HirYam}.  Just as in the Fuchsian case, Hiroe shows
that the existence of an irreducible framable connection with the given
formal types is equivalent to the existence of an irreducible quiver
representation with certain properties.

The quiver associated to a collection
$\bftype=(\ftype_1,\dots,\ftype_k)$ of unramified types is rather
complicated.  The set of vertices is partitioned into ``base points'' and
``path points'': $I=I_{\mathrm{base}}\sqcup I_{\mathrm{path}}$.  The path
points arise from the residue terms appearing in the formal types as
in \eqref{ltresform}.  The base points come from the purely irregular
part of (some of) the formal types.

As a warm-up, we describe the quiver $Q_{\mathrm{base}}^\ftype$ for a
single purely irregular formal type $\ftype$ as in \eqref{ltresform}.
The vertices are labeled by the simultaneous eigenspaces of $\ftype$,
so $I_{\mathrm{base}}^\ftype=\{1,\dots,\ell\}$.  There are only arrows
from $j$ to $j'$ if $j<j'$, in which case we draw
$\deg_{z^{-1}}(q_j(z)-q_{j'}(z))-1$ arrows.

\begin{example}
  If $\ftype=\left(\diag(a,b,b)z^{-3}+\diag(0,c,d)z^{-2}\right)\dzz$
  with $a,b,c\ne 0$, $a\ne b$, and $c\ne d$, then the corresponding
  quiver is
  
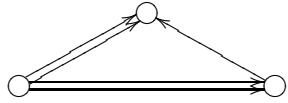
\[
\begin{xy}
(-17,0)*\cir<4pt>{}="A",
(17,0)*\cir<4pt>{}="B",
(0,9.7)*\cir<4pt>{}="C",
\ar@{->}@<.5mm> "A";"B",
\ar@{->}@<-.5mm> "A";"B",
\ar@{->}@<.5mm> "A";"C",
\ar@{->}@<-.5mm> "A";"C",
\ar@{->}"B";"C"
\end{xy}
\]
\end{example}

The following theorem explains the significance of this quiver in
studying moduli spaces of connections.

\begin{thm}[\cite{HirYam}]  Let $\ftype$ be a purely
    irregular formal type, and let $V=(V_1,\dots,V_\ell)$ be the
    eigenspace decomposition in \eqref{ltresform}.  Then there is a
    symplectic isomorphism between $\Ad^*(1+z\fgl(\pow))(\ftype)$ and
    $\Rep_{\overline{Q}_{\mathrm{base}}}(V)$.
  \end{thm}
  
  We now return to constructing $Q$ for
  $\bftype=(\ftype_1,\dots,\ftype_k)$.  In what follows, we call
  the index $i$ irregular (resp. regular) if $\ftype_i$ is, and we
  assume that $0$ is irregular.  Each $\ftype_i$ has a block diagonal
  decomposition, say with $\ell_i$ blocks, as in \eqref{ltresform}.
  We now let $I_{\mathrm{base}}=\{[i,j]|\text{$i=0$ or $\ell_i\ge
    2$, $1\le j\le \ell_i$}\}$ be the disjoint union of the
  $I_{\mathrm{base}}^{\ftype_i}$'s where either $i=0$ or $i$ is
  irregular with at least $2$ blocks.  We include all arrows in
  $Q_{\mathrm{base}}^{\ftype_i}$ and also draw arrows $[0,j]\to[i,j']$
  for every $i\in I_{\mathrm{base}}\setminus\{0\}$, $1\le j\le
  \ell_0$, $1\le j'\le \ell_i$.

  Each block in the decomposition of $\ftype_i$ gives rise to a
  residue term, so we obtain a collection of residues $R_{ij}$.  As in
  the Fuchsian case, each residue $R_{ij}$ has an associated path.  If
  $d_{ij}$ is the degree of its minimal polynomial, then this path
  consists of the $d_{ij}-1$ vertices $\{[i,j,k]\}$ with arrows
  $[i,j,k]\to[i,j,k-1]$.  We thus set $I_{\mathrm{path}}=\{[i,j,k]\}$.
  In the Fuchsian case, the path endpoints are all connected to the
  single base point $0$.  Here, the $[i,j,1]$'s will be connected to one
  or more vertices in $I_{\mathrm{base}}$.  If either $i=0$ or
  $\ell_i\ge 2$, then there is an arrow $[i,j,1]\to[i,j]$; otherwise,
  there are arrows $[i,1,1]\to[0,j]$ for $1\le j\le \ell_0$.

  The dimension vector $\alpha$ is defined just as in the Fuchsian
  case in terms of the Jordan block structure of the $R_{ij}$'s.
  Accordingly, $\alpha_{[i,j]}$ is the size the $j$th block for
  $\ftype_i$ while the $\alpha_{[i,j,k]}$'s are defined in terms of
  partial products of the minimal polynomial of $R_{ij}$. The
  deformation vector $\lambda$ is again defined in terms of the
  eigenvalues $\eta^k_{ij}$ of the $R_{ij}$.  For path vertices, it is
  just the same: $\lambda_{[i,j,k]}=\eta^k_{ij}-\eta^{k+1}_{ij}$.  For base
  vertices, we set $\lambda_{[i,j]}=-\eta^1_{ij}$ if $i\ne 0$ while
  $\lambda_{[0,j]}=-\eta^1_{0j}-\sum_{i\ne 0, \ell_i=1}\eta^1_{i1}$.

  We are finally ready to state Hiroe's theorem.  Let $L\subset\Z^I$
  be the sublattice defined by
  \begin{equation} L=\{\beta\in\Z^I\mid
    \sum_{j=1}^{\ell_0}\beta_{[0,j]}=\sum_{j=1}^{\ell_i}\beta_{[i,j]}\text{
      if $\ell_i\ge 2$}\}.
    \end{equation}  Recall that $p(\beta)=1-(1/2)\beta^t
    C_Q\beta$.
    
\begin{thm}[\cite{Hiroe}]  There is an irreducible framable connection
  with formal types $\bftype$ at the singular points
  $a_1,\dots,a_k$ if and only if
  \begin{enumerate} \item $\alpha$ is a positive root of $Q$, and
    $\alpha\cdot\lambda=0$;
    \item if $\alpha=\sum_{j=1}^\ell \beta_j$ with $\ell>2$ and
      $\beta_j\in L\cap  \Z_{\ge 0}^I$, then
      $p(\alpha)>\sum_{j=1}^\ell p(\beta_j)$.
    \end{enumerate}
  \end{thm}
  
In the Fuchsian case, Crawley-Boevey actually shows that there is an
isomorphism between $\sM^{\mathrm{st}}$ and the moduli space
$\sM^{\mathrm{reg}}_\lambda(Q,\alpha)$ of
irreducible representations of $Q$ with dimension vector $\alpha$ and
satisfying the deformation conditions \eqref{defcond}.  This is not
true in the general unramified case.  Indeed, Boalch has constructed
an example where the moduli space of irreducible framable connections
is not isomorphic to any quiver variety~\cite{Boa08}.  However, Hiroe
show that the stable moduli space is always isomorphic to an open
subset of a quiver variety.

\begin{thm}[\cite{Hiroe}]  There exists a vector
  $\lambda'\in\C^{I}$, agreeing with $\lambda$ on components in
  $I_{\mathrm{path}}$, for which there is an open embedding $\sM^{\mathrm{st}}\hookrightarrow\sM^{\mathrm{reg}}_{\lambda'}(Q,\alpha)$.  If
  there is at most one irregular unramified type, then $\lambda'=\lambda$
  and the open embedding is an isomorphism.
\end{thm}

\begin{rmk} Crawley-Boevey proved that the open embedding is an
  isomorphism when all formal types are regular singular~\cite{CB03}.
  The analogue for one irregular formal type is due to Hiroe and
  Yamakawa~\cite{HirYam}.
\end{rmk}

\section{Connections with ramified singular points}\label{s:ramsing}

\subsection{The Airy and Frenkel--Gross connections}\label{s:AFG}

In order to take ramified singularities into account, we need a notion
of \emph{ramified formal type}.  The LT normal form is not suitable as
it is not defined over the ground field $F$.  Instead, we need a
\emph{rational canonical form} for formal connections.

Recall that if $\ftype$ is an unramified formal type of slope $r$,
then $\ftype$ is an element of $\cB_r$, the principal parts of
bounded degree $r$.  The space $\cB_r$ is defined using the natural
degree filtration on the loop algebra and indeed can be identified
with $(z^{-r}\fgl_n(\pow)/z\fgl_n(\pow))\,\dzz$.  The polynomial loop
group $\GL_n(\pow)$ acts on $\cB_r$, and one then considers the
truncated orbit of $\ftype$ in $\cB_r$.

In order to discuss our basic problems in the ramified case, we need
analogues of these objects.  Thus, we need to associate the following
data to an arbitrary singularity:  
\begin{itemize}\item a formal type $\ftyperam$ giving the rational canonical form $d+\ftyperam$;
\item a subgroup $\Para\subset\GL_n(\pow)$ with Lie algebra $\fpara\subset\fgl_n(\pow)$;
\item a filtration $(\fpara^i)_{i\in\Z}$ of $\fgl_n(F)$  with
  $\fpara^0=\fpara$ which is
  preserved by conjugation by $\Para$;
\item a level $r$  in the filtration for which $\ftyperam$ may be viewed as an
  element of $\fpara^{-r}/\fpara^1\,\dzz$; and
\item the orbit of $\ftyperam$ in $\fpara^{-r}/\fpara^1\,\dzz$ under
  truncated conjugation by $P$ (or equivalently, the $P$-coadjoint orbit).
\end{itemize}

We begin by describing two important motivating examples.  First,
consider the classical Airy equation $y''=xy$.  This is the simplest
equation which has a turning point, namely a point where the solutions
(the \emph{Airy functions}) change from being oscillatory to
exponential.  For example, the solution of Schr\"odinger's equation
for a particle confined to a triangular potential well involves Airy
functions, reflecting the fact that the wave function of the particle
oscillates within the well, but then exhibits quantum tunneling into
the wall as it decays exponentially.  More generally, one can consider
the order $n$ version of the Airy equation: $y^{(n)}=xy$.  The Airy
equation is nonsingular except at $\infty$.  It will be convenient to
make the change of variables $z=x^{-1}$.  After converting the
equation into matrix form, one sees that the corresponding connection
is given by

\begin{equation}\label{airy}
\nabla_{\mathrm{Airy}}=d+z^{-1} 
\begin{pmatrix}
0 &  \cdots&0 & z^{-1}\\
1 & 0 & \cdots & 0 \\
\vdots & \ddots & \ddots & \vdots  \\
0  & \cdots & 1& 0 
\end{pmatrix}\dzz.
\end{equation}

Although the matrix of this connection has a pole of order $3$, the
matrix in the
leading term $z^{-2}e_{1n}\dzz$ is nilpotent, so one cannot conclude
that the slope at $z=0$  is $2$.  In fact, the
LT normal form is
$d+z^{-\frac{n+1}{n}}\diag(1,\zeta,\dots,\zeta^{n-1})\dzz$, where
$\zeta=e^{2\pi i/n}$, and the slope is the common value $\frac{n+1}{n}$
of all the slopes.  Thus, $\nabla_{\mathrm{Airy}}$ is ramified; in
fact, since the slope has
 largest possible denominator (in lowest form), we say that it is
 \emph{maximally ramified}.

 The Airy connection suggests how to construct a meromorphic connection with a
 maximally ramified singularity of slope $1/n$.  Let $\om_n=\sum_{i=1}^{n-1}e_{i,i+1}+ze_{n,1}$, i.e., $\om_n$ is the
matrix with $1$'s along the superdiagonal, $z$ in the
lower-left entry, and $0$'s elsewhere. Note that $\om_n^n=z$.  With
this notation, 
$\nabla_{\mathrm{Airy}}=d+z^{-1}\om_n^{-1}\dzz$.  We can now define
the rank $n$ Frenkel--Gross (or Kloosterman) connection by
\begin{equation}\label{fg} \nabla_{\mathrm{FG}}=d+\om_n^{-1}\sdfrac{dz}{z}.
\end{equation}
The corresponding differential equation is hypergeometric and was extensively studied by
Katz~\cite{KatzKloosterman,KatzExponential} (see also
Deligne~\cite{Del77}); it is a geometric
incarnation of a certain exponential sum called a Kloosterman sum
which arises in number theory.  Unlike the Airy connection, it has two
singular points, one regular singular and the other irregular.  The irregular
singular point has slope $1/n$ and LT normal form
$d+z^{-1/n}\diag(1,\zeta,\dots,\zeta^{n-1})\dzz$ while the
regular singular point has principal unipotent monodromy.  As we will
see later, both the Frenkel--Gross and Airy connections are
\emph{rigid}, i.e., uniquely determined by their local behavior.

For these two connections, we will view the matrix of the connection
in \eqref{airy} and \eqref{fg} as the ramified formal type of the
irregular singularity.  The group associated to these formal types
is the ``standard Iwahori subgroup'' $\Iwa$, consisting of those
$g\in\GL_n(\pow)$ which are upper triangular modulo $z$.  The
corresponding Iwahori filtration $\{\iwa^s\}$ will be defined below,
but it will satisfy $\om_n^s\in\iwa^s\setminus\iwa^{s+1}$.  This
means that $\om^{-1}\dzz$ may be viewed as the leading term of the
formal Frenkel--Gross connection $\fc_{\mathrm{FG}}$ with respect to
this filtration, and similarly for the formal Airy connection.
The Iwahori 
filtration is periodic with period $n$ in the sense that
$\iwa^{s+n}=z\iwa^s$, and as we will see, these facts show that
the irregular singularities of $\nabla_{\mathrm{Airy}}$ and
$\nabla_{\mathrm{FG}}$ have slopes $\frac{n+1}{n}$ and $1/n$
respectively, without needing to compute the LT forms.  

The matrix $\om_n$ is regular semisimple (though not diagonalizable
without extending scalars), and this means that the leading terms of
$\fc_{\mathrm{Airy}}$ and $\fc_{\mathrm{FG}}$ are regular semisimple
as well.  Formal connections with ``regular semisimple leading term''
with respect to an appropriate filtration are called \emph{toral}.  To
avoid technicalities, we focus here on toral formal types. However, a
rational canonical form for arbitrary ramified connections may be
obtained from Sabbah's refined Levelt--Turrittin
decomposition~\cite{Sab08}; we will return to this in a future paper.

\subsection{Lattice chain filtrations and parahoric subgroups}\label{s:lattice}

As we have seen, ramified formal types should be associated to certain
filtrations of the loop algebra.  For simplicity, we will restrict
attention to filtrations arising from ``lattice chains''; this will be
enough to discuss toral formal types~\cite{BrSa1,BrSa5}. It
will also suffice to define a notion of the ``leading term'' of a
formal connection when the naive leading term with respect to the
usual valuation is not useful.

An \emph{$\pow$-lattice} $L$ in $F^n$ is a finitely generated
$\pow$-module satisfying $L\otimes_{\pow}\laur \cong \laur^n$.
Concretely, it is the $\pow$-span of some basis of $F^n$. A
\emph{lattice chain} in $\laur^n$ is a collection $(L^i)_{i\in\Z}$ of
lattices satisfying the following properties:
\begin{enumerate}
\item $L^i\supsetneq L^{i+1}$ for all $i$; and
\item there exists a positive integer $e$, called the \emph{period}, such that $L^{i+e} = zL^i$ for all $i$.
\end{enumerate}
A \emph{parahoric subgroup} $\Para\subset\GL_n(\laur)$ is the
stabilizer of a lattice chain, i.e.,
$\Para=\{g\in\GL_n(\laur) \mid gL^i = L^i\text{ for all }i\}$.  We
let $e_{\Para}$ denote the period of the lattice chain stabilized
by $\Para$.  

We will be particular interested in parahoric subgroups associated to
lattice chains with nice uniformity properties.
\begin{definition} A lattice chain  $(L^i)_{i\in\Z}$ (and the
  corresponding parahoric subgroup) is called
\begin{enumerate} \item \emph{uniform} if $e|n$ and $\dim L^i/L^{i-1}=n/e$
for all $i$;
\item \emph{almost uniform} if $e|(n-1)$, $e\ne 1$, $\dim
L^i/L^{i-1}=1+(n-1)/e$ for $i$ in a single residue class modulo $e$,
and $\dim L^i/L^{i-1}=(n-1)/e$ otherwise;
\item \emph{complete} if its period is as large
as possible, i.e., $e=n$.
\end{enumerate}
\end{definition}
Note that a complete lattice chain is necessarily uniform.  The parahoric subgroups associated to complete lattice chains are
called \emph{Iwahori subgroups}.

The group $\GL_n(F)$ acts on the set of lattice chains with finitely
many orbits parametrized by subsets of $\Z/n\Z$ modulo the translation
action.  To see this, first note that any lattice $L$ has a valuation $\nu(L)\in\Z$: if
$L$ has an $\pow$-basis $f_1,\dots,f_n$, then the valuation of the
determinant of the matrix whose columns are the $f_i$'s is independent
of the choice of basis.  Since $\nu(zL)=z^n\nu(L)$, it is clear that a
lattice chain consists of lattices with valuations determined by some
subset  $J\subset \Z/n\Z$.  If we let $X_J$ be the set of all such
lattice chains, then $\SL_n(F)$ acts transitively on
$X_J$. (These are the (partial) affine flag varieties for
$\SL_n$~\cite{Sa00}.)  On the other hand, the action of $g\in\GL_n(F)$
will change the type of the lattice chain by translating $J$ by
$\nu(\det(g))$.

If $0\in J$, one can choose a lattice chain in $X_J$ whose
corresponding parahoric subgroup is block-upper triangular modulo $z$.
Let $0=k_0<k_1<\dots< k_{e-1}$ be the representatives of $J$ in
$[0,n-1]$; we set $k_e=n$.  Let $L^0=\pow^n$, and let $L^j$ be the
lattice with basis $\{ze_i\mid i>n-k_j\}\cup\{e_i\mid i\le n-k_j\}$
for $j>0$.  These lattices generate a lattice chain $\cL_J\in X_J$; we
denote its stabilizer by $P_J$ and call it the standard parabolic
subgroup of type $J$.  Let $Q_J\subset\GL_n(\C)$ be the block-upper
triangular parabolic subgroup with blocks of sizes $k_j-k_{j-1}$
arranged in order of decreasing $j$.  Then $P_J=\pi^{-1}(Q_J)$, where
$\pi:\GL_n(\pow)\to\GL_n(\C)$ is the homomorphism induced by the
``evaluation at $0$'' map $z\mapsto 0$.

Uniform standard parahoric subgroups correspond to subgroups of
$\Z/n\Z$; all of their blocks have the same size.  If $e|n$, we let
$P_e$ denote the uniform standard parahoric subgroup of period $e$.  In
particular, the trivial subgroup corresponds to the maximal parahoric
subgroup $P_1=\GL_n(\pow)$ while $J=\Z/n\Z$ gives the standard Iwahori
subgroup $P_n=I$: the pullback of the standard Borel subgroup of upper
triangular matrices.  If $e|(n-1)$, but $e\ne 1$, we will similarly
let $P_e$ denote the almost uniform standard parabolic subgroup with
$J$ given by the residues of $\{j\frac{n-1}{e}\mid 0\le j<e\}$.

Each lattice chain $\{L^i\}_{i\in\Z}$---say with corresponding
parahoric subgroup $\Para$---determines a filtration $\{\fpara^s\}_{s\in\Z}$
of $\fgl_n(\laur)$:
\begin{equation}\fpara^s = \{X\in\fgl_n(\laur) \mid XL^i \subset L^{i+s}\text{ for
    all }i\};
\end{equation}
in particular, $\fpara^0=\fpara\coloneqq\Lie(\Para)$.  One also gets a
filtration $\{\Para^s\}_{s\in\Z_{\ge 0}}$ of $\Para$ defined by
$\Para^0=\Para$ and $\Para^s=1+\fpara^s$ for $s>0$.  In the special
case that each lattice in the lattice chain stabilized by $\Para$
admits an $\pow$-basis of the form $\{z^{k_j}e_j\}_{j=1}^n$, the
corresponding filtration $\{\fpara^s\}_i$ is induced by a grading
$\fgl_n(\C[z,z^{-1}])=\bigoplus_s \fpara(s)$.  In particular, this is
true for the standard lattice chains $\cL_J$.

When $\Para$ is uniform, this filtration has a particular simple
description.  Indeed, there exists an element $\om_P\in\GL_n(F)$ such
that $\fpara^s=\om_P^s\fpara=\fpara \om_P^s$~\cite{Bu}.  In the
standard case $P=P_J$, one can choose $\om_J=\om_{P_J}$ so that it
generates the grading as well.  This means that if $\ell_J$ is the
usual block-diagonal Levi subalgebra of $\fq_J=\Lie(Q_J)$, then
$\fpara(s)=\ell_J\om_J^s=\om_J^s\ell_J$.  In particular, if $I$ is the
standard Iwahori subgroup, then $\om=\om_n$ from the previous section
has the desired property: $\iwa^s=\om^s\iwa=\iwa\om^s$ and
$\iwa(s)=\om^s \ft=\ft\om^s$, where $\ft$ is the diagonal Cartan
subalgebra.  In the same way that $\om$ only has nonzero entries on
the ``first circulant superdiagonal'', in the general standard uniform
case, the element $\om_P$ will have nonzero entries on the first
circulant block superdiagonal.

\begin{example}
  Let $n=3$.  Then
  \[\om_3 = \begin{pmatrix} 0 & 1 & 0 \\ 0 & 0 & 1 \\ z & 0 & 0 \end{pmatrix}\and\Iwa=\begin{pmatrix} \pow^* & \pow & \pow \\ z\pow & \pow^* & \pow \\ z\pow & z\pow & \pow^* \end{pmatrix}.\]
  Some steps in the standard Iwahori filtration on $\fgl_n(\laur)$ are shown below:
  \[
    \iwa^{-2}=\begin{pmatrix} \pow & z^{-1}\pow & z^{-1}\pow \\ \pow &
      \pow & z^{-1}\pow \\ \pow & \pow &
      \pow \end{pmatrix}\quad\subsetneq\quad \iwa^{-1}=\begin{pmatrix}
      \pow & \pow & z^{-1}\pow \\ \pow & \pow & \pow \\ z\pow & \pow &
      \pow \end{pmatrix}\quad\subsetneq\quad \iwa^{0}=\begin{pmatrix}
      \pow & \pow & \pow \\ z\pow & \pow & \pow \\ z\pow & z\pow &
      \pow \end{pmatrix}.\]
\end{example}

\subsection{Fundamental strata}\label{FS}
  
Let $\fc=d+(M_{-r}z^{-r}+M_{-r+1}z^{-r+1}+\dots)\dzz$ be a formal
connection with $r>0$.  If $M_{-r}$ is non-nilpotent, then the slope
of $\fc=r$.  However, if the slope of $\fc$ is not an integer, then
$M_{-r}$ will be nilpotent regardless of the trivialization for the
underlying vector bundle.  In this case, the native leading term
$M_{-r}z^{-r}\dzz$ of the connection is not useful.

A better notion of the ``leading term'' of $\fc$ was introduced by
Bremer and Sage in \cite{BrSa1} in terms of a geometric theory
of fundamental strata, motivated by the analogous $p$-adic theory of
Bushnell and Kutzko~\cite{Bu,BuKu}.  A \emph{$\GL_n$-stratum} is a
triple $(P,s,\beta)$ consisting of a parahoric subgroup, a nonnegative
integer $s$, and a functional $\beta\in(\fp^s/\fp^{s+1})^\vee\cong
\fp^{-s}/\fp^{-s+1}\,\dzz$.  We say that the stratum has depth
$s/e_P$. 
An element $\tilde{\beta}\in \fp^{-s}$ inducing $\beta$ is
called a \emph{representative} of $\beta$.  If the filtration
associated to $P$ comes from a grading (for example, if $P$ is a
standard parahoric subgroup), then there is a unique homogeneous
representative $\beta^\flat\in\fp(-s)$.  The stratum is called
\emph{fundamental} if every representative is nonnilpotent or
equivalently, if $\beta^\flat$ is nonnilpotent (when it is defined).

Given a formal connection $(\hV,\fc)$, we let $[\fc]_\phi$ denote the
matrix of the connection with respect to the trivialization $\phi$.
In  other words, 
$\fc=d+[\fc]_\phi$ in this gauge.   We say that $\fc$ contains a
positive depth stratum with respect to $\phi$ if
$[\fc]_\phi/\dzz\in\fp^{-s}$ and is a representative of $\beta$.  (See
\cite{BrSa1,BrSa3} for the general definition.)  Intuitively, this
means that $\beta$ is the leading term of the connection with respect
to the filtration $\{\fpara^j\}$.

\begin{thm}[\cite{BrSa1}, Theorem 4.10]\label{fsgl} Any formal connection
  $(\hV,\fc)$ contains a fundamental stratum
  $(P,r,\beta)$ with $P$ a standard parahoric subgroup.  The
  depth is positive if and only if $\fc$ is irregular singular, and
  $\slope(\fc)=r/e_P$.  Moreover, 
\begin{enumerate} \item If $\fc$ contains the stratum
  $(P',r',\beta')$, then $r'/e_{P'}\ge \slope(\fc)$.
\item If $\fc$ is irregular singular, then a stratum $(P',r',\beta')$
  contained in $\fc$ is fundamental if and only if $r'/e_{P'}= \slope(\fc)$.
\end{enumerate}
\end{thm}

For example, an unramified connection of slope $r$ in Levelt--Turrittin
form \eqref{LT} contains the fundamental stratum
$(\GL_n(\pow),r,D_rz^{-r}\dzz)$.  The Frenkel--Gross and Airy formal
connections contain the fundamental strata $(I,1,\om^{-1}\dzz)$ and
$(I,1,\om^{-(n+1)}\dzz)$ respectively.  This gives another proof that
their slopes are $1/n$ and $(n+1)/n$.

\subsection{Regular strata and toral connections}

Classically, much more information is available about a formal
connection when its (naive) leading term is regular semisimple and not
just nonnilpotent.  The analogue for strata involves the notion of a
\emph{regular stratum}.  A stratum $(P,r,\beta)$ is called regular if
every representative $\tilde{\beta}$ is regular semisimple.  If this
is the case, the centralizers of the various representatives are
maximal tori in the loop group which are conjugate under $P^1$.  We
abuse notation and call the stratum $S$-regular if $S$ is one of these
maximal tori.  If the homogeneous representative $\beta^\flat$ is
defined, then the stratum is regular if and only if $\beta^\flat$ is
regular semisimple.  We call a formal connection \emph{$S$-toral} if
it contains an $S$-regular stratum.  For example, an unramified formal
connection is $T(F)$-toral (with $T$ the standard diagonal maximal
torus) if and only if the leading term of the LT normal form is
regular semisimple.

Most classes of maximal tori do not give rise to $S$-toral
connections.  Recall that there is a bijection between the conjugacy
classes of maximal tori in $\GL_n(F)$ and conjugacy classes in the
symmetric group $\symmetric{n}$, hence to partitions of
$n$~\cite{KaLu88}.  If $\mu=(\mu_1,\dots,\mu_k)$ is a partition, then
the corresponding maximal tori are those isomorphic to $\prod_{i=1}^k
F[z^{1/\mu_i}]^\times$.  Such a maximal torus may be realized
explicitly as the block-diagonal embedding $\prod_{i=1}^k
\C(\!(\om_{\mu_i})\!)^\times$, where we set $\om_1=z$.  We call a
maximal torus \emph{uniform} if each block has the same size $e$ (so
$n=ek$). It is \emph{almost uniform} if it splits as $S'\times
F^\times$ where $S'$ is a uniform maximal torus in $\GL_{n-1}(F)$.  It
turns out that there is an $S$-regular stratum $(P,r,\beta)$ when
\begin{enumerate*}\item $S$ is uniform with $e$ blocks and
  $\gcd(r,e)=1$; or \item $S=S'\times F^\times$, $S'$ is a uniform
  maximal torus in $\GL_{n-1}(F)$ with $e$ blocks, and $\gcd(r,e)=1$.
\end{enumerate*}
Moreover, $P$ can be taken conjugate to a standard
parahoric subgroup whose block structure corresponds to the partition
dual to $\mu$.  (For example, one can take $P$ uniform with period $e$
in the first case.)  In particular, we conclude that toral connections
necessarily have slope $r/e$ with $\gcd(r,e)=1$ and $e$ dividing $n$
or $n-1$.

We will be particularly interested in the case of maximally ramified
toral connections, i.e., toral connections of slope $r/n$ with
$\gcd(r,n)=1$.  There is a unique class of maximal tori in
$\GL_n(\laur)$ that are anisotropic modulo the center, meaning that
they have no non-central rational cocharacters.  Concretely, such tori
are as far from being split as possible.  They are called
\emph{Coxeter maximal tori}, since they correspond to the Coxeter
class in $\symmetric{n}$ consisting of the $n$-cycles.  A specific
representative of this class is the standard Coxeter torus
$\cC=\C(\!(\om_n)\!)^\times$ with Lie algebra $\fcox$.  The formal
Frenkel--Gross and Airy connections are both $\cC$-toral.

In order to discuss toral formal types, we now introduce
representatives of the classes of uniform and almost uniform Cartan
subalgebras that are well-behaved with respect to the parahoric
subgroups $P_e$ defined in section~\ref{s:lattice}.  If $\fs$ is a
uniform Cartan subalgebra isomorphic to $F[z^{1/e}]^{n/e}$, it is
endowed with a filtration $\fs^j$ obtained by assigning degree one to
$z^{1/e}$.  (There is also an associated grading with graded pieces
$\fs(j)$.)  We say that $\fs$ is \emph{compatible} to $\fpara$ if
$\fs^j=\fpara^j\cap\fs$ for all $j$.  In the case that $\fpara$ is
standard, we say that $\fs$ is \emph{graded compatible} to $\fpara$ if
$\fs(j)=\fpara(j)\cap\fs$ for all $j$.

The block-diagonal uniform Cartan subalgebra with blocks of
size $e$ is not compatible with $\fpara_e$.  Instead, set
$V_j=\spa\{e_i|i\equiv j\pmod{e}\}$, so
$F^n=V_1\oplus\dots\oplus V_{n/e}$. Let $\om_{ej}$ be our usual matrix
$\om_e$ with respect to the given graded basis on $V_j$.  We now set
$\fs_e=\bigoplus_{j=1}^{n/e}\C(\!(\om_{ej})\!)$; it is graded compatible
with $\fp_e$.  Note that if $\gcd(r,e)=1$, then an element of
$\fs_e(r)$ is regular semisimple if it is of the form
$(c_j\om_{ej}^{r})_j$, where the $c_j$'s are nonzero constants such
that $c_j/c_{j'}$ is not an $e$th root of $1$ for $j\ne j'$.

If $e|(n-1)$, one can similar define $\fs_e$ graded
compatible with $\fp_e$; it is just the direct sum of a uniform
Cartan subalgebra in $\fgl_{n-1}(F)$ with $Fe_{nn}$.

The crucial feature of toral connections is that they can be
``diagonalized'' into a Cartan subalgebra.  If $(r,e)=1$, then the
space of $S_e$-formal types of slope $r/e$ is the subset
$\cA(S_e,r)\subset \left(\fs(-r)+\dots+\fs(0)\right)\dzz$ consisting
of matrices with regular semisimple component in degree $-r$.  It may
also be identified with a Zariski-open subset of $\fs^{-r}/\fs^1\dzz$.
\begin{thm}{\cite{BrSa1},\cite{BrSa5}}  If $\fc=d+[\fc]\dzz$ contains
  an $S_e$-regular stratum $(S_e,r,\beta)$, then $[\fc]\dzz$ is
  $P_e^1$-gauge equivalent to a unique element of $\cA(S_e,r)$.
\end{thm}
Thus, we get a rational canonical form for toral connections.

When $e=1$, so $S_e=T(F)$, this is just the familiar fact that if
$\nabla=d+(D_{-r}z^{-r}+M)\dzz$ with $D_{-r}\in\fgl_n(\C)$ regular
semisimple and $M\in z^{1-r}\fgl_n(\pow)$, then it is gauge equivalent
to a diagonal unramified LT normal form.  At the opposite extreme
when $e=n$, then $S_1=\cC$ and \begin{equation}
  \cA(\cC,r)=\{p(\om_n^{-1})\mid p\in\C[x], \deg p=r\}.
\end{equation}
In fact, one can show that all formal connections of slope $r/n$ with $\gcd(r,n)=1$ are
$\cC$-toral, so we have obtained a rational canonical form for
maximally ramified connections.
\begin{thm}[{\cite{KLMNS1}}]\label{diaggl}  Let $\fc$ be a maximally ramified connection
    of slope $r/n$ with $\gcd(r,n)=1$.  Then $\fc$ is formally gauge
    equivalent to a connection of the form $d+p(\om^{-1})\dzz$ with
    $p$ a polynomial of degree $r$.
  \end{thm}

  \begin{rmk}\label{waff} Let $S_e^0=\prod_{j=1}^{n/e}\C[\![\om_{ej}]\!]^\times$;
    it is the unique maximal bounded subgroup of $S_e$.  Let
    $W_{S_e}^{\mathrm{aff}}=N(S_e)/S_e^0$ be the relatively affine
    Weyl group of $S_e$.  It is the semidirect product of the relative
    Weyl group $W_{S_e}$ and the free abelian group $S_e/S_e^0$.
    There is a natural action of $W_{S_e}^{\mathrm{aff}}$ on
    $\cA(S_e,r)$: $W_{S_e}$ acts by conjugation in the usual way while
    $S_e/S_e^0$ acts by translations on $\fs_e(0)$.  It can be shown
    that the moduli space of $S_e$-toral connections of slope $r/e$ is
    given by $\cA(S_e,r)/W_{S_e}^{\mathrm{aff}}$~\cite{BrSa1}.
  \end{rmk}
  
We can now defined the truncated orbit of a toral formal type.  Fix
$\ftyperam\in\cA(S_e,r)$.  By
compatibility of $\fp_e$ and $\fs_e$, there are natural injections
$\fs_e^{-r}/\fs_e^1\,\dzz\hookrightarrow\fp_e^{-r}/\fp_e^1\,\dzz\subset\fp_e^\vee$.  The
parahoric subgroup
$P_e$ thus acts on $\fp_e^{-r}/\fp_e^1\,\dzz$ by the coadjoint action, and we
define the truncated orbit of $\ftyperam$ to be
$\orb_\ftyperam=\Ad^*(P_e)(\ftyperam)$.  Equivalently, one can view $\ftyperam$ as
an element of $\fp_e^{-r}$, and then $\orb_\ftyperam$ is obtained by
conjugation by $P_e$ followed by truncation.  Note that this action factors
through the finite-dimensional group $P_e/P_e^{r+1}$.

\subsection{Moduli spaces of toral framable connections}

We are now ready to discuss moduli spaces of meromorphic connections
with toral and unramified formal types.  First, we need to understand
what ``framability'' means in the context of toral singularities.
\begin{definition}\label{D:framable} Let $\gc=d+[\gc]_\phi$ be a global connection on $\pp$ with a singular point at $0$,
and let $\ftyperam$ be an $\Stor{e}$-formal type of depth $r$. We say
that $\gc$ is \emph{framable at $0$ with respect to $\ftyperam$} if \begin{enumerate}\item  there
  exists $g\in \GL_n(\C)$ such that $[\gc]_{g\phi}=g\cdot [\gc]_\phi\in\fpara_{e}^{-r}\dzz$ and
  $[\gc]-\ftyperam\in\fpara_{e}^{1-r}\dzz$; and
\item there exists an element $p\in \Para_{e}^1$ such that $\Ad^*(p)([\gc]_{g\phi}\dzz)|_{\fpara_e}=\ftyperam$.
\end{enumerate}
\end{definition}
Intuitively, this means that there is a global trivialization in which
the leading terms of the matrix of $\gc$ and $\ftyperam$ agree and that
there is a further formal gauge change by an element of $P_e^1$ in
which the matrix of $\gc$ and $\ftyperam$ are equal, viewed as functionals
on $\fp_e$.

Fix three disjoint (possibly empty) subsets $\{a_1,\dots,a_m\}$,
$\{b_1,\dots,b_\ell\}$, and $\{c_1,\dots,c_p\}$ of $\pp$ with $m\ge
1$.  Let $\bftyperam=(\ftyperam_1,\dots,\ftyperam_m)$ be a collection of
toral formal types at the $a_i$'s, $\bftype=(\ftype_1,\dots,\ftype_\ell)$
a collection of nontoral unramified irregular formal types at the
$b_j$'s, and $\borb=(\orb_1,\dots,\orb_p)$ a collection of nonresonant
adjoint orbits at the $c_k$'s.  One can now consider the category
$\sC(\bftyperam,\bftype,\borb)$ of meromorphic connections $\gc$ satisfying the
following properties:
\begin{enumerate}
\item $\gc$ has irregular singularities at the $a_i$'s and $b_j$'s, regular
  singularities at the $c_k$'s, and no other singular points;
\item for each $i$ and $j$, $\gc$ is framable at $a_i$ (resp. $b_j$)
  with respect to the formal type $\ftyperam_i$ (resp. $\ftype_j$); and
\item  for each $k$, $\gc$ has residue at $c_k$ in $\orb_k$.
\end{enumerate}

As in the Fuchsian and unramified cases, one can construct the
corresponding moduli space $\M(\bftyperam,\bftype,\borb)$ as the Hamiltonian
reduction of a product of symplectic $\GL_n(\C)$-manifolds, one for
each singular point.  It remains to define the symplectic manifold
$\M_\ftyperam$ encoding the local data at a toral singularity.  It no
longer suffices to take the truncated orbit
$\orb_\ftyperam=\Ad^*(P_e)(\ftyperam)$ because it no longer contains
full information about the principal part.  Indeed, $\fs_e(i)\dzz$
contributes to the residue for positive $i$ as long as $i<e$.  

Recall that $P_e$ is the pullback of a certain standard
parabolic subgroup $Q_e\subset\GL_n(\C)$ under the map
$\GL_n(\pow)\to\GL_n(\C)$ induced by $z\mapsto 0$.  More
intrinsically, $Q_e\cong P_e/(1+\GL_n(\pow))$.  The ``extended orbit''
$\M_\ftyperam\subset
(Q_e\backslash\GL_n(\C))\times\fgl_n(\pow)^\vee$ is now defined by
\begin{equation}\M_{\ftyperam}=\{(Q_e g,\alpha) \mid
  (\Ad^*(g)\alpha)|_{\fpara_e}\in \Ad^*(\Para_e)(\ftyperam)\}.
\end{equation}
Thus, $\alpha$ is the principal part, and the global trivialization $g$
is needed to move the principal part into $\fp_e^{-r}\dzz$.  One only
considers $Q_e$-orbits of such $g$ since $Q_e$ preserves
$\fp_e^{-r}\dzz$. 
It is easy to see that $\M_{\ftyperam}$ is a symplectic manifold~\cite{BrSa1}.  The group $\GL_n(\C)$ acts on $\M_{\ftyperam}$ via
$h\cdot(Q^b g,\alpha)=(Q^b gh^{-1},\Ad^*(h)\alpha)$; this action is
Hamiltonian with moment map
$(Q^b g, \alpha)\mapsto\alpha|_{\fgl_n(\C)}=\Res(\alpha)\dzz$.

\begin{thm}\label{genmodspace} The
  moduli space $\M(\bftyperam,\bftype,\borb)$ is given
  by 
\begin{equation*}  \M(\bftyperam,\bftype,\borb)\cong \left[ \left( \prod_i
      \M_{\ftyperam_i}\right) \times \left( \prod_j
      \orb_{\ftype_j}\right)\times \left( \prod_k \orb_k \right)
  \right] \sslash_0 \GL_n(\C).
\end{equation*}
\end{thm}

This theorem was proved by Bremer and Sage in \cite{BrSa1} with
the hypothesis that all irregular singular points are uniform toral.
The proof in this general case follows by combining the methods of
\cite{BrSa1} and \cite{HirYam}.  Note that as in the previous cases,
there is no dependence on the specific points.

\begin{rmk} When all irregular singular points are toral, there are
  two other variants of this moduli space.  One can consider the
  moduli space $\widetilde{\M}(\bftyperam,\borb)$ of \emph{framed
    connections}, which classifies framable connections together with
  the additional data of ``compatible framings''.  This was
  constructed by Boalch in the unramified case~\cite{Boa} and by
  Bremer and Sage when ramified uniform toral singularities are
  allowed~\cite{BrSa1}.  This moduli space is a smooth symplectic
  manifold, and $\M(\bftyperam,\borb)$ is obtained from it by
  Hamiltonian reduction by a certain torus.  There is also a larger
  Poisson moduli space classifying framed connections ``with fixed
  combinatorics''.  This means that at each $a_i$, the formal type is
  allowed to be arbitrary in $\cA(S_{e_i},r_i)$ for fixed $e_i$ and
  $r_i$.  Again, the construction is due to Boalch in the unramified
  case~\cite{Boa} and to Bremer and Sage in the uniform toral
  case~\cite{BrSa2}.
\end{rmk}

As usual, one can also define the stable moduli space
$\M^{\mathrm{st}}(\bfF,\bfD,\bfO)$ consisting of the irreducible
framable connections with the given formal types.  The \emph{toral
  Deligne--Simpson problem} is then the determination of when the stable
moduli space is nonempty.

We now compute these moduli spaces in the simplest case that involves a
ramified formal type: rank $2$ framable meromorphic connections
with a Coxeter toral singularity of slope $1/2$ and (possibly) a
regular singular point.  As we will see in the next section,
meromorphic connections with a Coxeter toral singular point are always
irreducible, so here, the stable moduli space coincides with the whole
moduli space.

\begin{example}\label{rameg}  Let $\ftyperam=a\om^{-1}_2\dzz$ with $a\in\C^*$, and let
  $\orb\subset\fgl_2(\C)$ be a nonresonant adjoint orbit.  The moduli
  space $\M(\ftyperam,\orb)$ is the space of $\GL_2(\C)$-orbits of
  triples $(Bg,\alpha,Y)$ with $(Bg,\alpha)\in\sM_\ftyperam$,
  $Y\in\orb$, and $\Res(\alpha)+Y=0$.  (Here, $B$ is the group of
  upper triangular matrices.)  This is the same as the space of
  $B$-orbits of triples $(B,\alpha,Y)$. 

  We first compute the truncated orbit of $\ftyperam$, using the fact
  that $I=T\rtimes I^1$ where $I^1 =1+\ft\om+\iwa^2$.  Since
  \begin{equation*}(1+\diag(b,c)\om)a\om^{-1}(1+\diag(b,c)\om)^{-1}=a\om^{-1}+\diag(b-c,c-b)\pmod{\iwa^1},
  \end{equation*}
  we get
\begin{equation}\label{fgorbit}
  \Ad^*(I)(\ftyperam)=\Ad^*(T)\left\{\begin{pmatrix} u & az^{-1}\\a & -u
  \end{pmatrix}\sdfrac{dz}{z}
  \mid u\in\C\right\}.
\end{equation}
Thus, in order to satisfy the condition
$\alpha\dzz|_{\iwa}\in\orb_{\ftyperam}$, we must have 
\begin{equation}\alpha=\begin{pmatrix} u& vz^{-1}+w\\v^{-1} & -u
  \end{pmatrix}\sdfrac{dz}{z}
\end{equation}
for some $u,v,w\in\C$ with $v\ne 0$.  In fact, each $B$-orbit has a
unique representative with $v=1$ and $u=0$.  This means that the only
adjoint orbits $\orb$ that give nonempty moduli space are those containing
$\left(\begin{smallmatrix} 0&-w\\-1&0
  \end{smallmatrix}\right)$ for some $w$.  Thus, $\M(\ftyperam,\orb)$ is a
singleton if $\orb$ is regular nilpotent or regular semisimple with
trace zero; otherwise, it is empty.  We remark that in the regular
nilpotent case, the unique such connection is the $\GL_2$-version of
the Frenkel--Gross rigid connection, and this argument shows that it is
indeed uniquely determined by its local behavior.
\end{example}

\subsection{The ramified Deligne--Simpson problem for Coxeter connections}\label{S:ramDS}

Little is known about the ramified Deligne--Simpson problem in general.
We now restrict attention to a simple special case in which the
problem has recently been solved by Kulkarni, Livesay, Matherne,
Nguyen, and Sage~\cite{KLMNS1,KLMNS2}: connections with a
\emph{maximally ramified} irregular singularity (i.e., an irregular
singularity of slope $r/n$ with $\gcd(r,n)=1$) and (possibly) an
additional regular singular point.  Without loss of generality, we can
view such connections as connections on $\Gm$ with the irregular
singularity at $0$.  Following \cite{KS2}, we refer to such
connections as \emph{Coxeter connections}.  As we have seen, classical
examples include the Frenkel--Gross and Airy connections.  Coxeter
connections have played a significant role in recent work on the
irregular geometric Langlands
program~\cite{FrGr,HNY13,Zhu17,KS2,JKY21}. They have also arisen in
Lam and Templier's proof of mirror symmetry for minuscule flag
varieties~\cite{LamTemp17}.

It is shown in Theorem 4.4 of \cite{KLMNS1} (c.f. \cite{KS3}) that
maximally ramified formal connections are precisely the $\cC$-toral
connections.  Moreover, $\cC$-toral connections are
irreducible~\cite{KS2}, so meromorphic connections with a $\cC$-toral
singularity are also irreducible. Thus, the Deligne--Simpson problem
for Coxeter connections is the determination of when the moduli space
$\M(\ftyperam,\orb)$ is nonempty for $\ftyperam$ a $\cC$-formal type
and $\orb$ a nonresonant adjoint orbit.
 
\begin{thm}{\cite[Theorem 5.4]{KLMNS1}}\label{DScox}  Let $\ftyperam$ be a maximally
ramified formal type of slope $r/n$ with $\gcd(r,n)=1$.  Then
$\M(\ftyperam,\orb)$ is nonempty if \begin{enumerate}\item
  $\Tr(\orb)=-\Res\Tr(\ftyperam)$, and
\item there are at most $r$ Jordan blocks for each eigenvalue of
  $\orb$.
\end{enumerate}
\end{thm}

We can reformulate this result in a more geometric way.  The set $\Orb^n$ of
adjoint orbits in $\fgl_n(\C)$ is a partially ordered set under the Zariski closure
relation:  $\orb\preceq\orb'$ if and only if
$\orb\subset\overline{\orb'}$.  The set of matrices with fixed
characteristic polynomial $q$ is a closed subset of $\fgl_n(\C)$ which
is stable under conjugation, so as posets $\Orb^n=\bigsqcup\orbp{q}^n$, where
the disjoint union runs over all monic polynomials of degree $n$.

Let $q=\prod_{i=1}^s(x-c_i)^{m_i}$ be a polynomial of degree $n$ with
the $c_i$'s distinct complex numbers.  Let $\prt(m)$ be the set of
partitions of $m$; it is a poset under the dominance order.  By the
theory of the Jordan canonical form, there is a bijection
$\orbp{q}^n\xrightarrow{\sim}\prod_{i=1}^s\prt(m_i)$ sending an orbit to the
partitions determined by the sizes of the Jordan blocks for each
eigenvalue.  Under this identification, the closure ordering corresponds to the product of the dominance orderings. 

We will show that the subset
\begin{equation}\DS(\ftyperam, q)=\{\orb\in\orbp{q}^n\mid \M(\ftyperam, \orb) \neq \varnothing
  \}
\end{equation}
is a principal filter in $\orbp{q}^n$.  (This means that it is an
upward-closed subset with a unique minimum element.)  Given positive
integers $r$ and $m$, there exists a unique smallest partition
$\prn{r}{m}\in\prt(m)$ with at most $r$ parts.  Indeed, if we write
$m=kr+r'$ with $k,r\in\Z$ and $0\le r'<r$, then this partition is
given by $\prn{r}{n}= \{(k+1)^{r^\prime}, k^{r-r^\prime}\}$.  We now
let $\prp{r}{q}$ be the orbit in $\orbp{q}^n$ corresponding to the
element
\[(\prn{r}{m_1}, \prn{r}{m_2}, \ldots, \prn{r}{m_s})\in\prod_{i=1}^s
  \prt(m_i).\]  This is the unique smallest orbit for which each
  eigenvalue has at most $r$ eigenvalues.

\begin{example}\label{homog} For nilpotent orbits (i.e., $q=x^n$), a
  representative of $\prp{r}{x^n}$ is given by the matrix with $1$'s
  on the $r$th subdiagonal if $1\le r\le n-1$ and $0$ if $r\ge n$.
  Note that these matrices are precisely the residues of $\om_n^{-r}\dzz$.
\end{example}

\begin{corollary}\label{DScox2} Let $\ftyperam=p(\om^{-1}_n)\dzz$ be a maximally ramified
  formal type of slope $r/n$, and let $q=\prod_{i=1}^s(x-c_i)^{m_i}$
  be a degree $n$ polynomial with the $c_i$'s distinct modulo $\Z$.
  Then,
\begin{enumerate}\item $\DS(\ftyperam, q)$ is empty if and only if
  $np(0)\neq -\sum_{i=1}^s m_i c_i$.
\item If $\DS(\ftyperam, q)$ is nonempty, it is the principal filter
  generated by $\prp{r}{q}$.
\item  For fixed $q$, $\DS(\ftyperam, q)$ only depends on $\deg p$ and
  $p(0)$.
\end{enumerate}
\end{corollary}

\begin{example}\hfill
\begin{enumerate}[label=(\alph*)] \item  If $\slope(\ftyperam)=1/n$, then
    $\M(\ftyperam,\orb)$ is empty unless the trace condition holds and
    $\orb$ is regular, i.e., maximal in $\Orb^n$.
\item If $\slope(\ftyperam)=(n-1)/n$, then $\M(\ftyperam,\orb)$ is
  nonempty unless the trace condition fails or $\orb=\{aI\}$ for some
  $a\in\C$.
\item If $\slope(\ftyperam)>1$, then  $\M(\ftyperam,\orb)$ is
  nonempty unless the trace condition fails.   More generally, this is
  the case if $\slope(\ftyperam)\ge m_i/n$ for all $i$, with
  $q=\prod_{i=1}^s(x-c_i)^{m_i}$ as above.
\end{enumerate}
\end{example}

\begin{rmk} In general, it is not easy to find explicit connections in
  $\M(\ftyperam,\orb)$.  However, it can be done if $\ftyperam$ is
  \emph{homogeneous}, meaning that $\ftyperam=a\om_n^{-r}\dzz$ and
  $\orb$ is nilpotent.  Indeed, an algorithm is described in
  \cite{KLMNS2} which constructs a binary, strictly upper triangular
  matrix $b_\orb$ for any $\orb\succeq\prp{r}{x^n}$ such that
  $d+(a\om_n^{-r}+b_\orb)\dzz$ gives an element of
  $\M(a\om_n^{-r}\dzz,\orb)$.  In particular, $b_{\prp{r}{x^n}}=0$.
\end{rmk}

\subsection{Rigidity for Coxeter connections}\label{rigidcoxgl}

We now discuss rigidity for Coxeter connections.  As in the case of
Fuchsian connections, there are several possible notions of rigidity
for connections with formal types as in Theorem~\ref{genmodspace}:
physical rigidity, ``moduli space rigidity'' (where the moduli space
is a singleton), and ``stable moduli space rigidity''.  In general, it
is very difficult to determine directly whether a connection is rigid
in any of these senses.  However, it turns out that there is a more
accessible version called \emph{cohomological rigidity}, which
coincides with physical rigidity for irreducible connections.

Let $j:U\hookrightarrow\pp$ be the inclusion of a nonempty
Zariski-open set.  A connection $\gc$ on $U$ is called cohomologically
rigid if $H^1(\pp,j_{!*}\ad_{\gc})=0$.  Here, $\ad_{\gc}$ is the
associated adjoint vector bundle.  For irreducible connections, this
cohomology vanishing condition implies that the connection has no
infinitesimal deformations.  It is a theorem of Bloch and Esnault that
cohomological and physical rigidity coincide for irreducible
connections~\cite{BE04}.  We will simply call such connections rigid.

We can now give a complete classification of rigid framable Coxeter
connections under the additional hypothesis that the connection has
unipotent monodromy at $\infty$.

\begin{thm}\label{rigidthm} Let $\ftyperam=p(\om_n^{-1})\dzz$ be a rank $n$ maximally ramified formal
  type of slope $r/n$ with $p(0)=0$, and let $\orb$ be any nilpotent
  orbit with at most $r$ Jordan blocks.  Then there exists a rigid
  connection with formal type $\ftyperam$ and unipotent monodromy
  determined by $\orb$ if and only if $\orb=\prp{r}{x^n}$ and $r|(n\pm
  1)$.
\end{thm}

In particular, the Frenkel--Gross and Airy connections are rigid. Note
that there are no rigid Coxeter connections of this type with the
slope of the irregular singularity greater than $(n+1)/n$.

\begin{rmk} It was originally proved in \cite{KS2} that
  $d+\om_n^{-r}\dzz$ is rigid if and only if $r|(n\pm 1)$.
\end{rmk}

We conclude this section by giving a brief sketch of the proof.  Let
$\gc$ be a Coxeter connection.  There is an explicit formula for the
cohomology group $\dim(H^1(\pp,j_{!*}\ad_{\gc}))$ in terms of the
global and local differential Galois groups of $\gc$ and of the
irregularity of the adjoint connection at $0$.  Let $\Galgp$ denote
the global differential Galois group, and let $\Galgp_0$ and
$\Galgp_\infty$ denote the local differential Galois groups at $0$ and
$\infty$.  These are all closed algebraic subgroups of $\GL_n(\C)$.
Also, let $\Irr(\ad_{\fc_0})$ denote the irregularity of the formal
connection $\ad_{\fc_0}$.  It is shown in \cite[Proposition 11]{FrGr}
that
\begin{equation}\label{numcrit} \dim(H^1(\pp,j_{!*}\ad_{\gc}))=\Irr(\ad_{\fc_0})-\dim(\fgl_n(\C)^{\Galgp_0})
-\dim(\fgl_n(\C)^{\Galgp_\infty})+2\dim(\fgl_n(\C)^{\Galgp}).
\end{equation}

It is a general fact that if $\fc$ is a toral connection of slope $s$,
then $\Irr(\ad_{\fc_0})=sn(n-1)$~\cite[Lemma 19]{KS1}.  Since $s=r/n$ in this
case, we have $\Irr(\ad_{\fc_0})=r(n-1)$.  

The differential Galois group of a $\cC$-toral connection is computed
in \cite{KS3}, following \cite{FrGr} and \cite{KS2}.  It has the form
$\Galgp_0\cong H\ltimes \langle\theta\rangle$, where $H$ is a
certain torus containing a regular semisimple element and $\theta$ is
an order $n$ element of $N(H)$.  One can then show as in
\cite{KS2} that $\fgl_n(\C)^{\Galgp_0}=\C$.  Since
$\Galgp_0\subset\Galgp$, we also get $\fgl_n(\C)^{\Galgp}=\C$.

Finally, if $N\in\orb$, then $\Galgp_\infty$ is the one-parameter
subgroup generated by $N$, i.e., the unipotent group $\{\exp(aN)\mid
a\in\C\}$.  It follows that $\fgl_n(\C)^{\Galgp_\infty}$ is the
centralizer of $N$.

A short calculation shows that \begin{equation}
  \dim(H^1(\pp,j_{!*}\ad_{\gc}))=\dim(\orb)+(r-n-1)(n-1)\ge 0.
\end{equation}
Recall that by Corollary~\ref{DScox2}, $\DS(\ftyperam,x^n)$ is the
principal filter generated by $\prp{r}{x^n}$, so in particular, there
exists a connection with formal types $\ftyperam$ and
$\orb$ if $\orb\succeq \prp{r}{x^n}$. If $\orb\succ \prp{r}{x^n}$, so
$\dim(\orb)>\dim(\prp{r}{x^n})$, then the dimension of the cohomology
group is positive.  Thus, rigidity is only possible if $\orb=
\prp{r}{x^n}$.  It remains to check when
$\dim(\prp{r}{x^n})=(n+1-r)(n-1)$.  It was shown in \cite{KS2} that
this occurs precisely when $r|(n\pm 1)$.

\section{Meromorphic $G$-connections}\label{s:gconn}

\subsection{$G$-connections}
A rank $n$ connection may be viewed as a
\emph{$\GL_n(\C)$-connection}: a principal $\GL_n(\C)$-bundle endowed
with a meromorphic connection.  This perspective allows us to
generalize our set-up to consider $G$-connections where $G$ is an
complex reductive group.  We can now ask the same fundamental
question: to what extent does specified local behavior determine a
global $G$-connection?  In this section, we briefly describe what is
needed to accomplish this as well as what is known.

Let $G$ be a connected reductive group over $\C$ with Lie algebra
$\fg$, and let $\Rep(G)$ be the category of finite-dimensional complex
representations of $G$.  Let $\cE$ be a principal $G$-bundle over
$\pp$.  By Tannakian formalism, $\cE$ may be viewed as a tensor
functor from $\Rep(G)$ to the category of (finite rank) vector bundles
via $V\mapsto V_\cE = \cE \times_G V$.  A \emph{meromorphic
  $G$-connection} is a tensor functor $(\cE,\gc)$ from $\Rep(G)$ to
the category of vector bundles endowed with a meromorphic connection.
We will always assume that $\cE$ is trivializable.  Fix a
trivialization $\phi$.  This induces corresponding trivializations
$\phi_V$ for each $V_\cE$.  The tensor functor is then determined by
$[\gc]_\phi\in \fg(\C(z))\dzz$; if $(V,\rho)$ is a representation,
then the induced connection on $V_{\cE}$ is given by $d+[\rho
([\gc]_\phi)]$.  We may thus view $\gc$ (in this trivialization) as
the abstract connection $d+[\gc]_\phi$.  Formal $G$-connections are
defined similarly.

Changing the trivialization changes the matrix of the $G$-connection
by gauge change:
\begin{equation}
  [\fc]_{g\phi}=g\cdot[\fc]_\phi=\Ad^*(g)([\fc]_\phi)-(dg)g^{-1}.
\end{equation}
Here, we are considering formal connections, and $g\in G(F)$; the
right invariant Mauer--Cartan form $(dg)g^{-1}$ lies in $\fg(F)$, so
the formula makes sense.  For meromorphic connections, one takes $g\in
G$, so the nonlinear gauge term is of course $0$.

We remark that there is an equivalence of categories between
$\GL_n$-connections and rank $n$ connections given by
$(\cE,\gc)\mapsto (V_{\cE},\gc_{V})$, where $V$ is the standard
representation of $\GL_n(\C)$.

For each $p\in\pp$, the local behavior of a meromorphic $G$-connection
$\gc$ at $p$ is governed by the corresponding formal connection
$\fc_p$.  Our original three questions about existence, uniqueness,
and moduli spaces of meromorphic connections with specified local
behavior now make sense for $G$-connections.  As we will see, for
Fuchsian and unramified $G$-connections, the construction of moduli
spaces may be generalized to arbitrary reductive $G$ with essentially
no changes.  However, very little is known about the analogues of the
Deligne--Simpson and rigidity problems.  It is more challenging to make
sense of the three problems in the ramified case, and we will describe
how to do so.  We will conclude the paper with a discussion of Coxeter
$G$-connections for $G$ simple.

Throughout this section, we fix a Borel subgroup and maximal torus
$G\supset B\supset T$, and let $\fb\supset\ft$ be the corresponding
Borel subalgebra and Cartan subalgebra.  Let $U$ be the unipotent
radical of $B$, and let $\fu$ be its Lie algebra.  We also fix a
nondegenerate invariant symmetric bilinear form $\left<\, , \right>$ on
$\fg$.  We use $\Res \left<\, , \right>$ to identify $\fg(F)\dzz$ with
the space of smooth $\C$-linear functionals $\fg(F)^\vee$.

\subsection{Fuchsian and unramified framable $G$-connections}

Just as for rank $n$ connections, in order to consider moduli spaces
of meromorphic connections with only regular singular points, we
restrict to \emph{Fuchsian $G$-connections}.  These are defined just
as before: One simply takes $G$-connections of the form
$\gc=d+\sum_{i=0}^k\frac{A_i}{z-a_i}dz$, with $A_i\in\fg$.  The
isomorphism classes of a regular singularity are parameterized by the
(local) monodromy conjugacy class. As long as each residue of the
Fuchsian connection $\gc$ satisfies an appropriate notion of
``nonresonance''\footnote{A formal $G$-connection $d+(B_0+B_1 z+\dots)\dzz$ is
  called \emph{nonresonant} if $\ad(B_0)\in\End(\fg)$ has no
  eigenvalues in $\Z\setminus\{0\}$.}, the local monodromies are determined by the
residue.  In other words, we can view
nonresonant adjoint orbits in $\fg$ as formal types for Fuchsian
connections.  We immediately get the analogue of
Proposition~\ref{p:fuchsmodspace}:
\begin{prop} Given nonresonant adjoint
orbits $\orb_1,\dots,\orb_k\subset\fg$, the corresponding moduli
space of Fuchsian $G$-connections is given by
\begin{equation}\sM(\orb_1,\dots\orb_k)=\{(A_1,\dots,A_k)\in\orb_1\times\dots\times\orb_k\mid\sum_{i=1}^k A_i=0\}/G.
\end{equation}
\end{prop}

The generalization of unramified framable connections to
$G$-connections is also straightforward.  There is a Jordan form for
formal $G$-connections due to Babbitt and Varadarajan that is entirely
analogous to the LT form for rank $n$ formal connections.

\begin{thm}[\cite{BV83}]  Any
  formal $G$-connection is isomorphic to a connection of the form
  \begin{equation}\label{BV} d+(D_rz^{-r/b}+D_{r-1}z^{-(r-1)/b}+\dots+
    D_1z^{-1}+D_0+N)\sdfrac{dz}{z},
\end{equation} where $D_i\in\ft$, $D_r\ne 0$, $N\in\fu$ commutes with each $D_i$, and $b\in\N$.
\end{thm}
The rational number $r/b$ is an invariant called the slope.

The formal connection $\fc$ is called unramified if one may take $b=1$
in the Jordan form, and we view the matrix of this Jordan form as an
unramified formal type.  Thus, the unramified formal types are subsets
of the space of principal parts $\fg(\pow)^\vee\cong
\fg(\C(z^{-1}))\dzz$.  The truncated orbits of an unramified formal type
$\ftype$ is just the coadjoint orbit $\Ad^*(G(\pow)(\ftype)\subset
\fg(\pow)^\vee$.

We again only consider \emph{framable $G$-connections}; if we assume
that $\infty$ is a regular point, these are connections of the form
\begin{equation}\label{eq:ppg} d+\bigg(\sum_{i=1}^k \sum_{\nu=0}^{r_i}
  \frac{A^{(i)}_{\nu}}{(z-a_i)^{\nu}}\bigg)\sdfrac{dz}{z}
\end{equation}
with $A^{(i)}_{\nu}\in\fg$ and $A^{(i)}_{r_i}$ nonnilpotent if
$r_i>0$.

We can construct the moduli space of framable connections with
unramified formal types $\ftype_i$ at the singular points $a_i$ just as
for $\GL_n(\C)$.  Indeed, one need only replace $\GL_n(\C)$ by $G$ in
\eqref{unrammod}.

We can now formulate the Fuchsian and unramified framable
Deligne--Simpson problem.  A  meromorphic connection is called
\emph{Lie-irreducible} if it does not admit a reduction of structure
to a proper parabolic subgroup.  For framable (and Fuchsian)
connections, this simply means that there is no proper parabolic
subalgebra of $\fg$ containing all of the coefficients
$A_\nu^{(i)}$ from \eqref{eq:ppg}.  The stable moduli
space is the subspace consisting of the Lie-irreducible connections with
the given formal types.

The Deligne--Simpson and unramified Deligne--Simpson problems for
$G$-connections are just the determination of when the relevant stable
moduli spaces are nonempty.  However, in spite of the obvious
parallelism with the situation for $\GL_n$,  I am not aware of any
significant results on these problems outside of type
$A$.

\subsection{Fundamental $G$-strata and toral $G$-connections}

As we have seen for $G=\GL_n(\C)$, if $\fc$ is a formal $G$-connection
of nonintegral slope, then the naive leading term of $\fc$ with
respect to any trivialization will be nilpotent.  One can generalize
the theory of fundamental strata discussed in \S\ref{FS} to
$G$-connections to obtain a more useful notion of the leading
term~\cite{BrSa3}.  It is also possible to define regular $G$-strata
and toral $G$-connections~\cite{BrSa5}.  In this context, lattice
chain filtrations are replaced by a class of Lie-theoretically defined
filtrations on $\fg(F)$ called \emph{Moy--Prasad filtrations}, which
are parameterized by points in the \emph{Bruhat--Tits building} $\cB$
for $G(F)$.

To minimize technicalities, from now on we assume that $G$ is simple
and simply connected.  The standard Iwahori subgroup $I\subset G(F)$
is the pullback of $B$ under the ``evaluation at $0$'' map $G(\pow)\to
G$.  A proper subgroup $P\subset G(F)$ is called a parahoric subgroup
if it contains a conjugate of $I$.  For example, parahoric subgroups
in $\SL_n(F)$ are stabilizers of lattice chains just as for $\GL_n$,
but here, the conjugacy classes of parahoric subgroups are determined
by subsets of $\Z/n\Z$ instead of classes of such subsets under the
translation action.  The Bruhat--Tits building is a simplicial complex whose
facets are in bijective correspondence with the parahoric subgroups of
$G(F)$.  The chambers (i.e., the maximal simplices) correspond to
Iwahori subgroups while vertices correspond to maximal parahoric
subgroups.  We call the closure of the chamber corresponding to $I$
the \emph{fundamental alcove}.  The building is divided into
\emph{apartments}--real affine spaces of dimension equal to the rank
of $G$--which are parameterized by the split maximal tori of $G(F)$.
The apartment labeled by $T(F)$ is called the standard apartment; it
contains the fundamental alcove.

Given $x\in \cB$, let $G(F)_x$ (resp. $\fg(F)_x$) denote the parahoric
subgroup (resp. subalgebra) corresponding to the facet containing $x$.
The Moy--Prasad filtration $\{\fg(F)_{x,r}\}_{r\in\R}$ associated to
$x$ is a decreasing $\R$-filtration of $\fg(F)$ by $\pow$-lattices.
It satisfies $\fg(F)_{x,0}=\fg(F)_x$ and is $1$-periodic, in the sense
that $\fg(F)_{x,r+1}=z\fg(F)_{x,r}$. If we set
$\fg(F)_{x,r+}=\bigcup_{s>r}\fg(F)_{x,s}$, then the set of $r$ where
$\fg(F)_{x,r}\ne\fg(F)_{x,r+}$ is discrete.  If $x$ is in the standard
apartment, then this filtration is induced by a grading
$\fg(\C[z,z^{-1}])=\bigoplus_{r\in\R} \fg(F)(r)$.  We remark that for
$\SL_n$ or $\GL_n$, a lattice chain filtration with stabilizer $P$ is
just a rescaling of the Moy--Prasad filtration at $x_P$, the barycenter
of the facet labeled by $P$~\cite{BrSa3}.  For example,
$\iwa^s=(\fgl_n(F))_{x_I,s/n}$ and $\iwa(s)=(\fgl_n(F))_{x_I}(s/n)$. 

A $G$-stratum is  a triple $(x,r,\beta)$ with $x\in\cB$,
$r\ge 0$, and \begin{equation*}\beta\in(\fg(F)_{x,r}/\fg(F)_{x,r+})^\vee\cong
\left(\fg(F)_{x,-r}/\fg(F)_{x,-r+}\right)\sdfrac{dz}{z}.
\end{equation*}
It is called
fundamental if every representative in this coset is nonnilpotent.
The definition of a formal $G$-connection containing a stratum is just
the same as for $\GL_n$.

We can now state the generalization of Theorem~\ref{fsgl}.

\begin{thm}[{\cite[Theorem 2.14]{BrSa3}}]
 Every formal $G$-connection $\fc$ contains a fundamental
  stratum $(x,r,\beta)$, where $x$ is in the fundamental alcove and $r\in\Q$; the
  depth $r$ is positive if and only if $\fc$ is irregular
  singular.  Moreover,
\begin{enumerate}
\item If $\fc$ contains the stratum $(y,r', \beta')$, then
$r' \ge r$.  
\item  If $\fc$ is irregular singular, a stratum
  $(y,r', \beta')$ contained in $\fc$ is fundamental if and only if
  $r' = r$.
\end{enumerate}
\end{thm}

The common depth of the fundamental strata contained in $\fc$ is one
possible definition of the slope of $\fc$.

To define toral $G$-connections, we need the notion of a regular
$G$-stratum.  For convenience, we will restrict to points in the
standard apartment, so that the functionals $\beta$ with always have a
homogeneous representative $\beta^\flat$.

Let $S\subset G(F)$ be a maximal torus.  The associated Cartan
subalgebra $\fs$ is endowed with a unique Moy--Prasad filtration
$\{\fs_r\}$ which comes from a grading with graded pieces $\fs(r)$.
We say that $x$ is compatible (resp. graded compatible) with $\fs$ if
$\fs_r=\fg_{x,r}\cap\,\fs$ (resp. $\fs(r)=\fg_{x}(r)\cap\fs$) for all
$r\in\R$.  The fundamental stratum $(x,r,\beta)$ with $x$ in the
standard apartment and $r>0$ is called \emph{$S$-regular} if $x$ is
compatible with $\fs$ and $S$ is the connected centralizer of
$\beta^\flat$.  If $\fc$ contains an $S$-regular stratum, then we say
that $\fc$ is \emph{$S$-toral}.

Just as for $\GL_n$, there is a bijection between conjugacy classes of
maximal tori in $G(F)$ and conjugacy classes in the Weyl
group~\cite{KaLu88}. A Weyl group element is called \emph{regular} if its
action on $\ft$ has a regular semisimple eigenvector; the
corresponding eigenvalue is also called regular.  If $\gamma$ is a
conjugacy class in $W$, we say that the maximal torus has type
$\gamma$ if it is in the corresponding class of tori.

\begin{thm}{\cite{BrSa5}}  Let $\gamma$ be  a conjugacy
  class in $W$.  Then the following are equivalent.
\begin{enumerate} \item $\gamma$ is regular and
  $e^{2\pi i r}$ is a regular eigenvalue for $\gamma$.
\item If $S$ is of type $\gamma$, then there exists a regular semisimple element in
$\fs(r)$.
\item There exists an $S$-toral connection with slope $r$ for some $S$
  of type $\gamma$.
\end{enumerate}
\end{thm}
For $\GL_n$, a maximal torus if of regular type if and only if it is
uniform or almost uniform.

Finally, we generalize Theorem~\ref{diaggl}. Let $\cA(S,r)$ be the open subset of
$\bigoplus_{j\in[-r,0]}\fs(j)$ whose leading component (i.e., the
component in $\fs(-r)$) is regular semisimple.  This is called the set
of \emph{$S$-formal types} of depth $r$.
\begin{thm}{\cite{BrSa5}}\label{formaltypes}
  If $\fc$ contains the $S$-regular stratum $(x,r,\beta)$, then $\fc$
  is $\fg_{x,0+}$-gauge equivalent to a connection with matrix in
  $\cA(S,r)\dzz$.
\end{thm}

Upon choosing nice representatives for the maximal tori admitting
toral connections, this theorem can be viewed as giving rational
canonical forms for formal connections.

\begin{rmk} As in Remark~\ref{waff}, there is an action of the relative affine
Weyl group $W^\mathrm{aff}_S$ on $\cA(S,r)$, and the moduli space of
$S$-toral connections of slope $e$ is given by
$\cA(S,r)/W^\mathrm{aff}_S$~\cite{BrSa5}.
\end{rmk}

\subsection{Moduli spaces of toral $G$-connections}

Let us return to our checklist of ingredients necessary to construct
moduli spaces.  We assume that $S$ is compatible with a point $x$ in
the standard apartment, and take $\ftyperam\in\cA(S,r)$. In fact, we
can (and will) choose $S$ so that it is graded compatible with a point
$x$ in the fundamental alcove.  The truncated orbit is then defined as
$\Ad^*(G(F)_x)(\ftyperam)$.  However, a difficulty now arises: the
parahoric subgroup $G(F)_x$ is not necessarily a subgroup of
$G(\pow)$.  This is problematic: principal parts are elements of
$\fg(\pow)^\vee$, and in the construction of moduli spaces, one needs
to be able to restrict these functionals to
$\fg(F)_x$ to check whether they give the desired formal type.  We thus
need to make the additional hypothesis that one can choose $x$ so that
the closure of the facet containing $x$ contains the vertex
corresponding to $\fg(\pow)$.  Note that this condition is always
satisfied for $\GL_n$, since $S_e$ is graded
compatible to the barycenter $x_{P_e}$.

If $x$ satisfies this hypothesis, then $G(F)_x$ is the pullback to $G(\pow)$ of a
standard parabolic subgroup $Q_x\subset G$.  We define the
extended orbit $\M_\ftyperam\subset
(Q_x\backslash G)\times\fg(\pow)^\vee$ just as before:
\begin{equation*}\M_{\ftyperam}=\{(Q_x g,\alpha) \mid
  (\Ad^*(g)\alpha)|_{\fg(F)_x}\in \Ad^*(G(F)_x)(\ftyperam)\}.
\end{equation*}
The construction of moduli spaces with toral singularities now goes
through as in Theorem~\ref{genmodspace}.  (See, for example,
\cite{Liv}, where this is worked out for $G=\Sp_{2n}(\C)$.)  The
ramified Deligne--Simpson problem for $G$-connections is the
determination of when the stable moduli space, consisting of the
Lie-irreducible connections with the given formal types, is nonempty.

\subsection{Coxeter $G$-connections}

Very little is known about the Deligne--Simpson problem for Fuchsian or
unramified $G$-connections.  Indeed, I do not even know of any
conjectural statements.  A major issue is that the conditions
appearing in the results of Crawley-Boevey and Hiroe are not
Lie-theoretic in nature.  Instead, they involve a translation to an
analogous problem for quiver representations, and this translation has
no known analogue outside of type $A$.  On the other hand, the
solution of the Deligne--Simpson (and rigidity) problem for Coxeter
connections is Lie-theoretic, involving the structure of the poset of
adjoint orbits, and thereby admits the possibility of a generalization
to other types.  Moreover, Coxeter $G$-connections can be defined for
any simple group (or for a reductive group with connected Dynkin
diagram)~\cite{KS2}.  In this final section, we discuss precise
conjectures (and some partial results) for the Deligne--Simpson and
rigidity problems for Coxeter $G$-connections.

Let $G$ be a simple group  with fixed Borel subgroup and maximal torus
$B\supset T$.  Let $\Phi$ denote the set of roots and $\Delta\subset
\Phi$ the set of simple roots.  We let $\fg_\alpha\subset\fg$ denote
the root space corresponding to the root $\alpha$; by convention, we
set $\fg_0=\ft$.  We fix an ``affine pinning'' of $\fg$, namely a
generator $X_\alpha\in\fg_\alpha$ for each negative simple root
$\alpha\in-\Delta$ and a generator $E\in\fg_\theta$, where $\theta$ is
the highest root.  We now set
$\om_{-1}=z^{-1}E+\sum_{\alpha\in-\Delta} X_\alpha$.  Note that for
$\SL_n$, an appropriate choice of the affine pinning gives
$\om_{-1}=\om_n^{-1}$ as defined in \S\ref{s:AFG}.

It can be shown that $\om_{-1}$ is regular semisimple~\cite{KS2}, and
we denote its connected centralizer by $\cC$.  We call $\cC$ the
Coxeter maximal torus and its Lie algebra $\fcox$ the Coxeter Cartan
subalgebra.  The terminology is justified by the fact that the
conjugacy class of $\cC$ corresponds to the Coxeter class in the Weyl
group $W$.

It turns out that $\cC$ is graded compatible with $x_I$, the
barycenter of the fundamental alcove.  Let $h$ be the Coxeter number
of $G$, i.e., the order of the Coxeter elements in $W$.  The
Moy--Prasad grading associated to $x_I$ has nonzero components only in
degrees $\frac{1}{h}\Z$.  It may be described explicitly in terms of
the heights of the roots of $\fg$:
\begin{equation}
  \fg_{x_I}(m+\ltfrac{s}{h})=z^m\left[\Bigg(\bigoplus_{\height(\alpha)=s}\fg_\alpha\Bigg)\oplus \Bigg(z\bigoplus_{\height(\alpha)=s-h}\fg_\alpha\Bigg)\right]\qquad\text{for $0\le s<h$.}
\end{equation}
The graded piece $\fcox(r/h)=\fg_{x_I}(r/h)\cap\fcox$ is
one-dimensional when $\gcd(r,h)=1$.  If $0<-r<h$, we
choose a generator $\om_r\in\fcox(r/h)$; in general,
$\om_{m+r/h}\overset{\text{def}}{=}z^m\om_r$ is a generator.
In type $A_{n-1}$, we simply set $\om_r=\om_n^r$.

If $r>0$ is relatively prime to $h$, then the set $\cA(\cC)$ of
$\cC$-formal types of depth $r/h$ is the subset of $\big(\bigoplus_{i=0}^r
\fcox(-i/h)\big)\dzz$ with nonzero component in degree $-r/h$. We call the
formal type \emph{homogeneous} if it has zero components outside of
this degree.

A meromorphic $G$-connection on $\Gm$ is called a \emph{Coxeter
$G$-connection} if it has a $\cC$-toral irregular singularity at $0$
and possibly an additional regular singular point at $\infty$.  If
$\ftyperam$ is a $\cC$-formal type and $\orb\subset\fg$ is a
nonresonant adjoint orbit, we can define the moduli space of framable Coxeter
$G$-connections $\M(\ftyperam,\orb)$  with these formal types.  Explicitly,
\begin{equation}\M(\ftyperam,\orb)\cong\{(\alpha,Y)\mid
\alpha\in\fg(\C[z^{-1}])\tdzz, Y\in\orb,
\alpha\vert_{\iwa}\in\Ad^*(I)(\ftyperam),\text{ and }\Res(\alpha)+Y=0\}/B.
\end{equation}
Coxeter connections are automatically Lie-irreducible because the
formal connection at $0$ is already Lie-irreducible.

We can now state our conjectural solution of the Deligne--Simpson
problem for Coxeter $G$-connections.  Let $\Orb^G$ be the poset of
adjoint orbits under the partial order determined by containment of
orbit closures.  Let $q$ be an element in the
adjoint quotient $\fg/G\cong \ft/W$. The elements in $\fg$ mapping to
$q$ form a $G$-stable closed subset; we let $\Orb^G_q$ be the collection
of adjoint
orbits contained in this subset. There is a poset isomorphism
$\Orb^G\cong\bigsqcup_q \Orb^G_q$. Let $\DS(\ftyperam,\orb)$ be the
set of  $\orb\in\Orb^G_q$ for which $\M(\ftyperam,\orb)$ is nonempty.

\begin{conjecture}\label{dscoxgen} Let $G$ be a simple complex group with Lie algebra
  $\fg$.  Fix a $\cC$-formal type $\ftyperam$ of slope $r/h$ with
  $\gcd(r,h)=1$.  Then for each $q\in\ft/W$, $\DS(\ftyperam,\orb)$ is
  either empty or is a principal filter generated by a (unique) orbit
  $\orb_q^r\in \Orb^G_q$ which is independent of the choice of
  $\ftyperam$ with the given slope.  Moreover,
\begin{enumerate} \item  If $r>h$, then $\orb_q^r=0$, so
  $\M(\ftyperam,\orb)$ is always nonempty.
\item If $q=0$ (so $\Orb^G_0$ is the set of nilpotent orbits), then
  $\Res(\om_{-r}\dzz)$ is a representative of $\orb_0^r$.
\end{enumerate}
\end{conjecture}

For $G$ simple, $\fcox(0)=0$; this explains the absence of the trace
condition which appears in Theorem~\ref{DScox}.

Finally, we discuss rigidity for Coxeter $G$-connections.  As in
\S\ref{rigidcoxgl}, we will actually investigate cohomological
rigidity, i.e., when $H^1(\pp,j_{!*}\ad_{\gc})=0$.  For
Lie-irreducible connections, the vanishing of this cohomology group
implies that $\gc$ has no infinitesimal deformations.  Unlike in type
$A$, it is not known whether physical rigidity and cohomological
rigidity coincide.  From now on, we will simply call cohomologically
rigid connections rigid.

We first consider \emph{homogeneous Coxeter connections}: meromorphic
$G$-connections of the form $\fc^a_r=d+a\om_{-r}\dzz$ with $r>0$
relatively prime to $h$ and $a\in\C^*$.  When $r=1$, we get the Frenkel--Gross rigid
$G$-connection.  In \cite{KS2}, Kamgarpour and Sage
showed that a homogeneous Coxeter $G$-connection with $r=h+1$ is
always rigid.  We call these connections \emph{Airy
  $G$-connections}.  Note that they are nonsingular at $\infty$.

It turns out that the Frenkel--Gross and Airy
connections are the only homogeneous Coxeter $G$-connections which are
rigid for all simple types.  The complete classification of rigid
homogeneous Coxeter $G$-connections is given in \cite{KS2}.

\begin{thm}\label{rigidcoxhom} The Coxeter
  connection $\nabla=\nabla_{r}^a$ is rigid if and only if
 \begin{enumerate} 
  \item[(i)] $r=h+1$ in which case
    $\nabla$ is regular at infinity;  
 \item[(ii)] $r=1$  in which case $\nabla$ is regular singular at infinity with principal unipotent monodromy; 
 \item[(iii)] $1<r<h$ satisfies the following conditions 
 \[
\begin{tabular}{|c|c|}
\hline
\textrm{Root system} & \textrm{Conditions on $r$} \\
\hline
$A_{n-1}$ & $r|n\pm 1$ \\
\hline
$B_n$ & $r|n+1, \,\, r|2n+1$\\
\hline
$C_n$ & $r|2n\pm 1$\\
\hline
$D_n$ & $r|2n, \, \, r|2n-1$\\
\hline
$E_7$ & $r=7$\\
\hline
\end{tabular}
\]
in which case $\nabla$ is regular singular at infinity with unipotent monodromy $\exp(\Res(\om_{-r}\dzz))$. 
 \end{enumerate} 
 \end{thm}  

We expect that the restriction to homogeneous formal types is
unnecessary.  Indeed, the following statement would follow from
the $q=0$ case of Conjecture~\ref{dscoxgen}.

\begin{conjecture}  Let $\ftyperam$ be a $\cC$-formal type of slope
  $r/h$, and let $\orb\subset\fg$ be any nilpotent orbit with
  $\orb\succeq\orb^r_0$ (as defined in the previous conjecture).
  Then, there exists a rigid connection with the given formal type and
  unipotent monodromy determined by $\orb$ if and only if $r$
  satisfies the conditions in Theorem~\ref{rigidcoxhom}.
\end{conjecture}


\bibliography{referencesconnlocbeh}
\bibliographystyle{amsalpha}

\end{document}